\documentclass[12pt]{amsart}
\usepackage [b5paper,left=1.5cm,right=1.5cm,top=2cm,bottom=2cm]{geometry}               
\usepackage{amssymb,amsmath,amsrefs}
\usepackage{enumerate}
\usepackage{color} 
\usepackage{graphicx}
\usepackage{mathrsfs}
\usepackage{setspace}
\usepackage{csquotes}  
\usepackage{hyperref}
\hypersetup{colorlinks=true,
            linkcolor=blue,
            anchorcolor=blue,
            citecolor=blue}

\usepackage{pifont}

\usepackage{caption}
\usepackage{subfigure}

\newtheorem{theo}{Theorem}

\newtheorem{lemma}{Lemma}

\newtheorem{prop}{Proposition}
\newtheorem{rema}{Remark}

\title{On the uniqueness of the planar 5-body central configuration with a trapezoidal convex hull}
\author[Yangshanshan Liu]{Yangshanshan Liu$^{1,2}$}
\author[Shiqing Zhang]{Shiqing Zhang$^{1}$$^\ast$}
\thanks{\noindent$^{1}$Chern Institute of Mathematics, Nankai University, Tianjin, China. 300071}
\thanks{\noindent$^{2}$Department of Mathematics, Sichuan University, Chengdu, China. 610065}
\thanks{$\ast$Corresponding author; E-mail: zhangshiqing@scu.edu.cn}
\date{}  
\begin{document}

\begin{abstract}
To apply Morse's critical point theory, we use mutual distances as coordinates to discuss a kind of central configuration of the planar Newtonian 5-body problem with a trapezoidal convex hull, i.e., four of the five bodies are located at the vertices of a trapezoid, and the fifth one is located on one of the parallel sides. We show that there is at most one central configuration of this geometrical shape for a given cyclic order of the five bodies along the convex hull. 
\end{abstract}

\maketitle

\section{Introduction}
The classical $n$-body problem is related with the study of $n$ particles $P_1,\cdots,P_n$ with masses $m_i>0$ and positions $x_i=(x_{i1},\cdots,x_{id})^T \in \mathbb{R}^{d}\,(i=1,\cdots,n;\, d=2,3)$ interacting with each other and satisfying the following differential equations
\begin{equation}\label{newton}
	m_i\ddot x_{ij} =\frac{\partial U}{\partial x_{ij}},\quad i=1,\cdots,n,\,j=1,\cdots,d,
\end{equation}
where 
\begin{equation*}
	U=\sum\limits_{i<j}\frac{m_im_j}{\vert x_j-x_i \vert}=\sum\limits_{i<j}\frac{m_im_j}{r_{ij}}
\end{equation*}
is the Newtonian potential and $r_{ij}=|x_i-x_j|$ is the mutual distance between $P_i$ and $P_j$. We denote by $x=(x_1,\cdots,x_n)\in\mathbb{R}^{dn}\backslash\Delta$ a collision free configuration, where $\Delta=\left\{x\in\mathbb{R}^{dn}|x_i=x_j,\forall \,i\neq j\right\}$ is the \textit{singular set}, and $r=(r_{12},r_{13},\cdots,r_{n-1,n})\in(\mathbb{R}^+)^{C_n^2}$ its corresponding mutual distance vector.

Central configurations are special arrangements of the $n$-body system satisfying a set of nonlinear algebraic equations
\begin{equation}\label{cc}
	-\lambda m_j(x_j-c)=\frac{\partial U(x)}{\partial x_j}=\sum_{i\neq j,i=1}^{n}\frac{m_im_j(x_i-x_j)}{r_{ij}^3},
\end{equation}
where $c=\frac{1}{m}\sum \limits_{i=1}^n m_ix_i$ is the \textit{center of mass} with $m=\sum\limits_{i=1}^nm_i$ the \textit{total mass} and $\lambda$ the \textit{configuration constant}. Furthermore, if we denote by 
\begin{equation}\label{moi}
I(x)=\frac{1}{2}\sum\limits_{i=1}^nm_i(x_i-c)^2=\frac{1}{2m}\sum\limits_{i<j}m_im_jr_{ij}^2
\end{equation}
the \textit{moment of inertia} of the system, a more compact form of \eqref{cc} is as follows
\begin{equation}\label{ccc}
	\nabla U(x)+\lambda\nabla I(x)=0,
\end{equation}
in which $\nabla$ denotes the gradient operator with respect to $x$.
Hence, $\lambda$ can be seen as the Lagrangian multiplier of the potential $U(x)$ with the constraint $I(x)=I_0$.
In short, if we simultaneously release the $n$ particles with zero velocity from this special position, all bodies will collide at the center of mass in a limited time. This can be seen as a special homothetic case of the \textit{homographic solutions} of \eqref{newton}, i.e.,
\begin{equation*}
	x(t)-c(t)=r(t)Q(t)(x_0-c_0),
\end{equation*}
where $x_0$ is a central configuration with $c_0$ its center of mass, $r(t)>0$ is a real scaling factor, namely, a solution of Kepler's problem (i.e., 2-body problem), and $Q(t)\in \mathrm{SO}(d)$ is a rotation. 
Furthermore, if we take $r(t)=1$, only considering the rotation, we get the well-known relative equilibrium solutions. 
These self-similar solutions may be the only analytic solutions to the $n$-body problems up to the recent research. The first two kinds of homographic solutions were constructed by Euler in 1767 and Lagrange in 1772 with initial particular positions for the 3-body problem, respectively, which came to be known as the Eulerian collinear and Lagrangian equilateral triangle center configurations. 

One can easily see that the equations \eqref{cc} and \eqref{ccc} are invariant under translation, dilation, and rotation, which allows us to consider the equivalent classes of central configurations naturally. In addition, the classical Morse theory requires the non-degeneracy of central configurations as critical points at least, which are not isolated under the continuous action of the rotation group. Hence, they are degenerate. 
Two central configurations are equivalent if they can be transformed from one to another by combining the above three transformations.
There are several open problems concerning central configurations, one of which is the finiteness conjecture, namely, the Chazy-Wintner-Smale conjecture \cites{chazy1918, smale1998, wintner1941}: \textbf{Is the number of equivalent classes of central configurations (or relative equilibria for $d=2$) finite for any given $n$ positive masses?} For $n=3$, Wintner \cite{wintner1941} showed that there are totally three Euler collinear central configurations and two Lagrange equilateral triangle central configurations. Hampton and Moeckel \cites{hampton2006} gave a positive answer for the case $n=4$, and it is still open for $n\geq5$ with several generic results on $n=5$, see \cites{albouy2012,hampton2011}. 
If we consider the special shapes of the configurations related to the classification of central configurations, some excellent results also come out. Moulton \cites{moulton1910} studied the $n$ collinear case in 1910, concluding that there are just $n!/2$ collinear central configurations for any given $n$ positive masses along a straight line. One can also refer to \cites{moeckel2015,smale1970b} for different proofs.
In other words, if the order of the $n$ masses along a line is fixed, the central configuration is unique.
Similarly, one can also get the uniqueness of some non-collinear central configurations provided the corresponding order of the masses is fixed; see \cites{cors2014,fernandes2020,xie2012} for some analytical approaches. Recently, Santoprete \cites{santoprete2021a,santoprete2021b} used a topological method to get the uniqueness results under the mutual distance coordinate system, different from the former ones, which can be summarized as follows 
\begin{enumerate}[(i)]
	\item\label{i} Use mutual distances as variables instead of positions to eliminate the rotation symmetry naturally;
	\item\label{ii} Replace the Cayley-Menger determinants constraint with simple new ones such that the gradients of both are paralleled; 
	\item\label{iii} Study the topology of a proper space formed by the corresponding constraints;
	\item\label{iv} Find out the critical point of the new Lagrangian function and compute the eigenvalues of the Hessian at that point;
	\item\label{v} Get the uniqueness of the non-degenerate critical point via Morse inequality.
\end{enumerate}
In this paper, we mainly apply this topological method to study a planar 5-body central configuration with a trapezoidal convex hull. Four of the five bodies are located on the vertices of a trapezoid, and the fifth is on one of the parallel sides. 
In 2012, Chen and Hsiao \cites{chen2012} discussed the existence of central configurations with this geometrical shape via some symmetric assumptions and numerical calculations. 
Here in this paper, we will show that for a fixed cyclic order of the five positive masses, there is at most one central configuration of this kind. 

It should be pointed out that (\ref{i}) and (\ref{ii}) are crucial for this approach. For (\ref{i}), the study of central configurations via mutual distances started from Dziobek \cites{dziobek1900}, who proved that the critical point equations 
\begin{equation}\label{dziobek4}
	\nabla_rU(r)+\lambda \nabla_rI(r)+\sigma \nabla_r\mathcal{F}_4(r)=0
\end{equation}
of the Lagrangian function $U(r)+\lambda (I(r)-I_0)+\sigma \mathcal{F}_4(r)$ with respect to $r$ are equivalent to \eqref{cc} for $n=4$, where
\begin{equation}\label{cayley-menger}
\mathcal{F}_n(1,2,\cdots,n)=\mathcal{F}_n(r)=
    \begin{vmatrix}0&1&1&\cdots&1\\
    1&0&r_{12}^2&\cdots&r_{1n}^2\\
    1&r_{12}^2&0&\cdots&r_{2n}^2\\
    \vdots&\vdots&\vdots&\ddots&\vdots\\
    1&r_{1n}^2&r_{2n}^2&\cdots&0
    \end{vmatrix}=(-1)^n\cdot 2^{n-1}\cdot ((n-1)!)^2\cdot \mathcal{V}_n^2,
\end{equation}
is the Cayley-Menger determinant and $\mathcal{V}_n$ is the corresponding volume of the simplex formed by $n$ points $P_1,\cdots,P_n$. One can refer to \cites{berger1987b} for more details. 
He also gave the famous Dziobek equations using the corresponding Cayley-Menger determinant as the geometrical constraint, in which
\begin{equation*}
	\frac{\partial \mathcal{F}_4(i,j,k,l)}{\partial r_{ij}^2}=-32\triangle_i\triangle_j,\quad i\neq j,
\end{equation*}
plays a key role (one can also refer to \cites{cors2012, hampton2002, moeckel2015} and \eqref{F4partialrij} in Subsection \ref{proof5co} in this paper),
where $\triangle_i$ denotes the oriental area of the triangle formed by the left three particles with the $i$-th deleted.
In comparison, the work of MacMillian and Bartky \cites{macmillan1932} on the planar 4-body problem, together with an extension work by Williams \cites{williams1938} on the planar 5-body problem, should be seen as a partial use of mutual distances, whose equations of central configurations are still based on the position coordinate system. Later in the 1980s, after the introduction of topological methods to the $n$-body problems by Smale \cites{smale1970b} and the further development of relative equilibria by Palmore \cites{palmore1973},
Meyer and Schmidt \cites{meyer1988b} as well as Schmidt \cites{schmidt1988, schmidt2002} found some new applications of the similar equations as \eqref{dziobek4}, or to say, the extrema equations of the Lagrangian function with Cayley-Menger determinants as constraints, and they got new central configurations by considering the perturbation of degenerate central configurations studied earlier by Palmore. 
For the cases of collinear 3-body, planar 4-body, and spacial 5-body, one Cayley-Menger determinant of $n$ particles $\mathcal{F}_n=0$ is enough to characterize them as Dziobek configurations, whose configuration dimension equals $n-2$.
But in \cites{schmidt2002}, Schmidt pointed out that for the planar 5-body problem, which is not a Dziobek one, any three Cayley-Menger determinants of four particles with one body deleted in total five choices are needed to guarantee that all the five bodies are coplanar, i.e., we need to consider the following equations
\begin{equation*}
	\nabla_r U(r)+\lambda \nabla_r I(r)+\sigma_1 \nabla_r \mathcal{F}_4(\hat i)+\sigma_2 \nabla_r \mathcal{F}_4(\hat j)+\sigma_3 \nabla_r \mathcal{F}_4(\hat k)=0,
\end{equation*} 
where $\mathcal{F}_4(\hat i)=\mathcal{F}_4(\hat j)=\mathcal{F}_4(\hat k)=0$ with $i\neq j\neq k$ and $\mathcal{F}_4(\hat i)$ denotes the Cayley-Menger determinant of four bodies with the $i$-th one deleted.
While Schmidt did not mention whether or not the above equations are equivalent to \eqref{cc} for $n=5$, and we try to give a proof, please see Proposition \ref{5co}.

For (\ref{ii}), replacing the Cayley-Menger determinant constraints with more simple ones can significantly simplify the calculations in (\ref{iii}), (\ref{iv}) and (\ref{v}),
and was first noticed by Cors and Roberts \cites{cors2012} in the study of co-circular central configurations of the 4-body problem in which they used the Ptolemy's theorem instead, supported by an important property, i.e., the gradients of both constraints with respect to $r$ are paralleled. Santoprete in \cites{santoprete2018} succeeded in replacing with a trapezoidal geometrical constraint, which can be obtained by the law of cosine of the four triangles separated by the diagonals of the trapezoid, and we will see in the following.

\section{The planar 5-body central configuration with a trapezoidal convex hull}

With the topological approach mentioned above, we can deal the cases without any symmetric assumptions.
We first show the equivalence of the equations of central configurations under both coordinate systems (Proposition \ref{5co}). Then, we claim that if it exists, this central configuration is a non-degenerate minimum of the potential function $U$ restricted to the corresponding manifold formed by the constraints (Lemma \ref{localmini}). Hence, we conclude that there is at most one central configuration of this kind:
\begin{theo}\label{5planarunique}
Given a cyclic (anti-clockwise) order of five positive masses, four of the five masses form a trapezoid, and the fifth is located on one of the parallel sides, not coinciding with the two vertices. Then, there is at most one central configuration for this type of geometric shape.
\end{theo}

\begin{figure}[htbp]
\begin{center}
\includegraphics[width=0.45\textwidth]{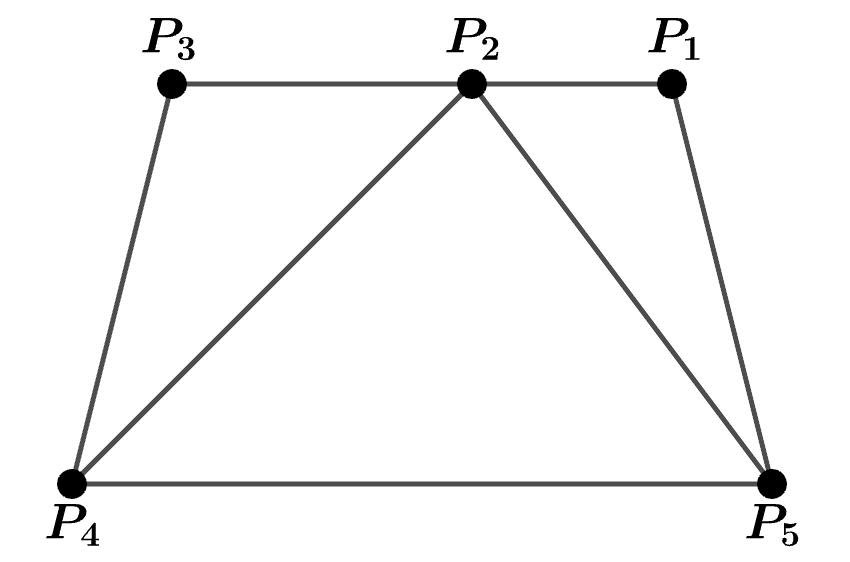}
\caption{A planar 5-body configuration with a trapezoidal convex hull.}
\label{trapezoid}
\end{center}
\end{figure}
\subsection{The reduced Lagrangian function}
Let 
$$x=(x_1,\cdots,x_5),x_i=(x_{i1},x_{i2})\in\mathbb{R}^2,i=1,\cdots,5$$ 
be the position vector and $r=(r_{12},r_{13},\cdots,r_{45})\in(\mathbb{R}^+)^{10}$ its distance vector, where
$$P_iP_j=r_{ij}=\vert x_j-x_i \vert=\sqrt{(x_{j1}-x_{i1})^2+(x_{j2}-x_{i2})^2},1\leq i<j\leq 5.$$ 
Without loss of generality, we fix an order of the five particles such that $P_1, P_2, P_3, P_4$ and $P_5$ lie anticlockwise along the trapezoidal convex hull. In addition, let $P_1P_3$ and $P_4P_5$ be the two parallel sides, with $P_2$ lying on $P_1P_3$, between the two endpoints; please see Figure \ref{trapezoid}.

According to \cites{schmidt2002} mentioned above, for the planar 5-body problem, any three Cayley-Menger determinants formed by four particles coming from five choices are needed to guarantee that all the five bodies are coplanar. Furthermore, with the observation from Cors and Roberts \cites{cors2012}, the introduction of \textit{geometrically realizable set} is necessary here too, i.e.,
\begin{equation}
\begin{aligned}
\mathcal{G}=\left\{r\in (\mathbb{R}^+)^{10}\vert \right.
&r_{ij}+r_{jk}\geq r_{ik}, \forall i,j,k \in\left\{1,2,3,4,5\right\},i\neq j\neq k,\\
&\left.\mathcal{F}_5(1,2,3,4,5)\leq0\right\},
\end{aligned}
\end{equation}
where $\mathcal{F}_5(1,2,3,4,5)\leq0$ can be obtained from \eqref{cayley-menger} with $n=5$.
Let 
\begin{equation*}
	\tilde W= U+\lambda (I-I_0)+ \eta_{i}F_i+\eta_{j}F_j+\eta_{k}F_k,\quad \forall i,j,k\in\{1,2,3,4,5\}
\end{equation*}
be the Lagrangian function, where we set $I_0=\frac{1}{2m}$ and
$F_i$ denotes the Cayley-Menger determinant formed by the four bodies with the $i$-th body deleted. Then we have 
\begin{prop}\label{5co}
 Suppose that $x=(x_1,\cdots,x_5)\in (\mathbb{R}^2)^5,x_i=(x_{i1},x_{i2})\in\mathbb{R}^2,i=1,\cdots,5$ is the position vector of the bodies in the planar 5-body problem and $r\in\mathcal{G}$ the corresponding distance vector. Then $x$ is the critical point of the potential $U$ restricted on the constraint $I=I_0$ if and only if $x$ is the critical point of $\tilde W(r(x))$ if and only if $r$ is the critical point of $\tilde W(r)$ provided there are no cases for 3-body collinear or 4-body collinear in the three chosen constraints $F_i$ in $\tilde W(r)$.
\end{prop}
We place the proof in Subsection \ref{proof5co}.

According to the position of the five bodies in Figure \ref{trapezoid}, we choose the constraints $F_1=0, F_2=0$, and $F_3=0$. 
Let
\begin{equation}
\begin{aligned}
	\mathcal{G}_2=\left\{\right.r\in(\mathbb{R}^+)^{10}\vert 
	& r_{ij}+r_{jk}> r_{ik},\forall\, i\neq j\neq k \in\left\{1,3,4,5\right\} \\
	 &\left.\mathcal{F}_4(1,3,4,5)\geq0
	\right\},
\end{aligned}
\end{equation}
which implies that any three particles of $P_1,P_3,P_4$ and $P_5$ are not collinear. 
Let
\begin{equation}
	\mathcal{N}_{123}=\left\{r\in\mathcal{G}\cap \mathcal{G}_2\vert I(r)=I_0,F_1=0,F_2=0,F_3=0\right\}
\end{equation}
be the normalized configuration space with respect to $r$. 
On the one hand, since points $P_2, P_3, P_4, P_5$, points $P_1, P_3, P_4, P_5$, and points $P_1, P_2, P_4, P_5$ are vertices of three trapezoids, we consider the geometrical constraints for each shape \cites{josefsson2013,santoprete2018}
\begin{subequations}
\begin{align}
T_1(r)&=2r_{23}r_{45}-r_{24}^2+r_{25}^2+r_{34}^2-r_{35}^2=0,\label{t1}\\
T_2(r)&=2r_{13}r_{45}-r_{14}^2+r_{15}^2+r_{34}^2-r_{35}^2=0,\label{t2}\\
T_3(r)&=2r_{12}r_{45}-r_{14}^2+r_{15}^2+r_{24}^2-r_{25}^2=0.\label{t3}
\end{align}
\end{subequations}
By substituting $r_{24}^2=2r_{23}r_{45}+r_{25}^2+r_{34}^2-r_{35}^2,r_{14}^2=2r_{13}r_{45}+r_{15}^2+r_{34}^2-r_{35}^2$ and $r_{14}^2=2r_{12}r_{45}+r_{15}^2+r_{24}^2-r_{25}^2$ into the Cayley-Menger determinants $F_1(r),F_2(r)$ and $F_3(r)$ respectively we obtain three remainders $-2K_1(r)^2,-2K_2(r)^2$ and $-2K_3(r)^2$ respectively, which implies with direct computation that 
\begin{equation}\label{substitrap}
 	F_i(r)=V_i(r)\cdot T_i(r)-2K_i(r)^2, \quad i=1,2,3,
 \end{equation}
 where 
\begin{equation*}
\begin{aligned}
V_1(r)=&-2[r_{35}^2(r_{23}-r_{45})^2-r_{24}^2r_{35}^2+(r_{23}^2-r_{25}^2)(r_{45}^2-r_{34}^2)],\\
K_1(r)=&r_{23}(r_{34}^2-r_{35}^2-r_{45}^2)+r_{45}(r_{23}^2-r_{25}^2+r_{35}^2),\\
V_2(r)=&-2[r_{35}^2(r_{13}-r_{45})^2-r_{14}^2r_{35}^2+(r_{13}^2-r_{15}^2)(r_{45}^2-r_{34}^2)],\\
K_2(r)=&r_{13}(r_{34}^2-r_{35}^2-r_{45}^2)+r_{45}(r_{13}^2-r_{15}^2+r_{35}^2),\\
V_3(r)=&-2[r_{25}^2(r_{12}-r_{45})^2-r_{14}^2r_{25}^2+(r_{12}^2-r_{15}^2)(r_{45}^2-r_{24}^2)],\\
K_3(r)=&r_{12}(r_{24}^2-r_{25}^2-r_{45}^2)+r_{45}(r_{12}^2-r_{15}^2+r_{25}^2).
\end{aligned}
\end{equation*}
Noticing that $P_1, P_2$ and $P_3$ are collinear, we denote the collinear constraint by
\begin{equation*}
	L_{1,2,3}(r)=r_{12}-r_{13}+r_{23}=0.
\end{equation*} 
Now let 
\begin{equation*}
\begin{aligned}
\mathcal{H}_{123}&=\left\{r\in\mathcal{G}\cap \mathcal{G}_2\vert T_1(r)=0,T_2(r)=0 ,T_3(r)=0 \right\},\\
\mathcal{H}_{T_2L}&=\left\{r\in\mathcal{G}\cap \mathcal{G}_2\vert T_2(r)=0 ,L_{1,2,3}(r)=0\right\}.	
\end{aligned}
\end{equation*}

We claim that $\mathcal{H}_{123}=\mathcal{H}_{T_2L}$. Actually, $\forall r\in \mathcal{H}_{123}$, by calculating we have $T_1+T_3-T_2=2(r_{12}+r_{23}-r_{13})r_{45}=0$, which implies $L_{1,2,3}=0$. Inversely, by Theorem 10 in \cite{josefsson2013}, for the convex quadrilateral formed by $P_2,P_3,P_4$ and $P_5$ we have $0=2r_{23}r_{45}\cos \xi-r_{24}^2+r_{25}^2+r_{34}^2-r_{35}^2\leq2r_{23}r_{45}-r_{24}^2+r_{25}^2+r_{34}^2-r_{35}^2=T_1(r)$, where $\xi$ denotes the angle between the extensions of the sides $r_{23}$ and $r_{45}$. Similarly we have $T_3(r)\geq0$. Therefore, $0=T_2=T_1+T_3$ forces the vanishing of both $T_1$ and $T_3$.

However, it is not easy to directly analyze the topological structure of $\mathcal{H}_{T_2L}$ since the existence of the geometrical realizable sets $\mathcal{G}$ and $\mathcal{G}_2$. To simplify, we consider some appropriate regions, as follows: 
We denote by 
\begin{equation*}
\begin{aligned}
	\mathcal{M}_{123}^+&=\left\{r\in(\mathbb{R}^+)^{10}\vert I(r)=I_0,T_1(r)=0, T_2(r)=0 ,T_{3}(r)=0\right\},\\
	\mathcal{M}_{T_iT_jL}^+&=\left\{r\in(\mathbb{R}^+)^{10}\vert I(r)=I_0,T_i(r)=0, T_j(r)=0 ,L_{1,2,3}(r)=0\right\}.
\end{aligned}
\end{equation*}
for $ i,j=1,2,3, i\neq j$
and 
\begin{equation*}
\mathcal{T}_{123}=\left\{r\in\mathcal{G}\cap\mathcal{G}_2\cap\mathcal{M}_{123}^+\vert I(r)=I_0,F_1(r)=0,F_2(r)=0,F_3(r)=0\right\}.
\end{equation*}
We claim that $\mathcal{M}_{123}^+=\mathcal{M}_{T_iT_jL}^+$ for any $i\neq j$ and $i,j=1,2,3$. Compared with $\mathcal{H}_{123}=\mathcal{H}_{T_2L}$ in the last paragraph, the different point here is the geometrical realizability among mutual distances. 
For arbitrary positive $r_{ij}$, the convex quadrilateral formula may not hold, i.e., they may not form an actual planar geometrical shape. Hence, the non-negativity of $T_i,i=1,2,3$ cannot be guaranteed. Therefore, one more study for $T_j$ when $j\neq i$ is needed. 
Without loss of generality, we take $i=1,j=3$, and the other two cases are similar. For any $r\in \mathcal{M}_{T_1T_3L}^+$, by calculating $T_1+T_3$ we have $2(r_{12}+r_{23})r_{45}-r_{14}^2+r_{15}^2+r_{34}^2-r_{35}^2=0$, and with $L_{1,2,3}=0$ we get $T_2=0$. The converse side is obvious with the equality $T_1+T_3-T_2=0$. 

The relationship of the above sets is shown in the following figure. The region we really want is $\mathcal{T}_{123}$, which is contained in the intersection of $\mathcal{H}_{123}$ and $\mathcal {M}_{123}^+$. In what follows we will show the topological structure of $\mathcal {M}_{123}^+=\mathcal{M}_{T_1T_3L}^+$ instead of $\mathcal{H}_{123}$.
\begin{figure}[htbp]
\begin{center}
\includegraphics[width=0.6\textwidth]{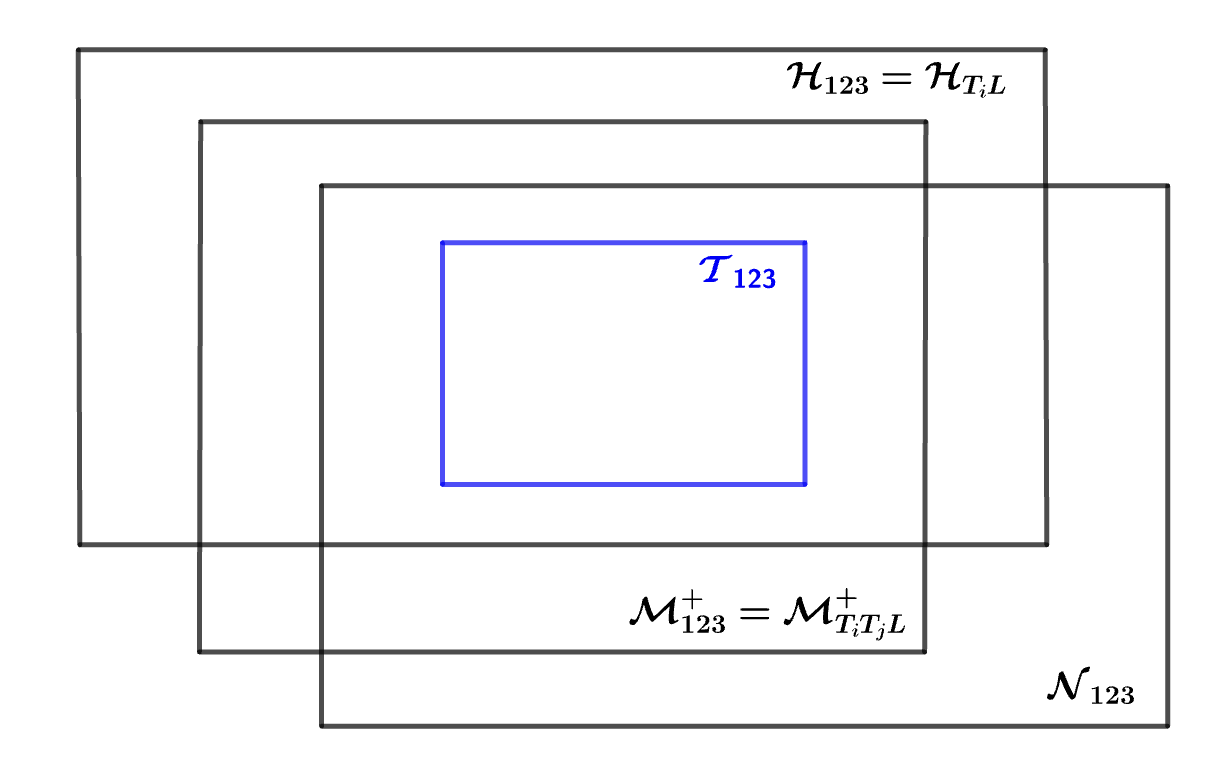}
\caption{$\mathcal{T}_{123}\subset \left(\mathcal{H}_{123}\cap \mathcal {M}_{123}^+\cap\mathcal{N}_{123} \right)$.}
\label{sets}
\end{center}
\end{figure}
\begin{lemma}\label{rF1F2F3}
If $r\in \mathcal{H}_{123}$, then $F_1=0,F_2=0$ and $F_3=0$.
\end{lemma}
\begin{proof}
$\forall r\in \mathcal{H}_{123}$, $\forall i\in \{1,2,3\}$ we have $T_i(r)=0$. From \eqref{substitrap} we have $F_i=-2K_i^2\leq 0$. From $r\in \mathcal{G}_2$ we have $F_i\geq 0$. Hence, we have $F_i=0$ for $i=1,2,3$.
\end{proof}
Now we can say that a configuration satisfying the position in Figure \ref{trapezoid} is a central configuration if and only if the corresponding distance vector $r\in\mathcal{T}_{123}$ is a critical point of the potential $U$ restricted on $\mathcal{N}_{123}$, i.e., $r\in\mathcal{T}_{123}$ is a critical point of 
\begin{equation*}
	\tilde W_{123}(r)=U(r)+\lambda (I(r)-I_0) +\eta_1 F_1(r)+\eta_2 F_2(r)+\eta_3 F_3(r).
\end{equation*}
In other words, $r\in\mathcal{T}_{123}$ satisfies 
\begin{equation}\label{nablatilde123}
	\nabla_r\tilde W_{123}(r)=\nabla_r U(r)+\lambda \nabla_r I(r)+\theta_1 \nabla_rF_1(r)+\theta_2 \nabla_rF_2(r)+\theta_3 \nabla_rF_3(r)=0.
\end{equation}
\begin{lemma}
For all $r\in\mathcal{H}_{123}$ and any $i\in \{1,2,3\}$, we have
\begin{equation}\label{allgradient}
	\nabla_r F_i(r)=V_i(r)\cdot \nabla_rT_i(r).
	\end{equation}
\end{lemma}
\begin{proof}
From \eqref{substitrap}, taking gradients with respect to $r$ we have
$$\nabla_r F_i =\nabla_rV_i\cdot T_i +V_i\cdot \nabla_r T_i-4K_i\cdot \nabla_r K_i.$$
Since $r\in\mathcal{H}_{123}$, we have $T_i=0$. Then by Lemma \ref{rF1F2F3} and \eqref{substitrap} we have $K_i=0$. This implies \eqref{allgradient}. 
\end{proof}
Substituting \eqref{allgradient} into \eqref{nablatilde123} we have 
\begin{equation*}\label{W123/r}
\nabla_r U(r)+\lambda \nabla_r I(r)+\eta_1\nabla_rT_1+\eta_2\nabla_rT_2+\eta_3\nabla_rT_3=0,
\end{equation*}
where
\begin{equation*}
		\eta_i=\theta_i\cdot V_i(r), \quad i=1,2,3.
\end{equation*}
This implies that $r\in\mathcal{T}_{123}$ can be seen as the critical point of the reduced Lagrangian function 
\begin{equation*}
	U(r)+\lambda(I(r)-I_0)+\eta_1T_1+\eta_2T_2+\eta_3T_3.
\end{equation*}
We have 
\begin{prop}\label{rinH123iffx}
$r\in\mathcal{T}_{123}$ is the critical point of $U$ restricted on $\mathcal{N}_{123}$ if and only if $r$ is the critical point of $U$ restricted on $\mathcal{M}_{123}^+$. Therefore, for a given cyclic order of the masses, $x$ is a central configuration satisfying the position assumption above if and only if its distance vector $r\in\mathcal{T}_{123}$ is a critical point of the Lagrangian function
	\begin{equation*}
	W_{123}(r)=U(r)+\lambda(I(r)-I_0)+\eta_1T_1(r)+\eta_2T_2(r)+\eta_3T_3(r)
\end{equation*}
where $I-I_0=0,T_i=0$ for $i=1,2,3$ with $\eta_i$ the Lagrangian multiplier.
\end{prop}
\begin{proof}
	From the above analysis, $\nabla_r\tilde W_{123}=0\Leftrightarrow \nabla_r W_{123}=0$ holds provided $r\in\mathcal{T}_{123}$.
\end{proof}
From $\nabla_r W_{123}=0$ and by direct computation we have the following ten equations
\begin{equation}\label{ten}
\left\{
\begin{matrix}
	S_{12}&&& &              &+& 2r_{45}\eta_3 &=0,\\
	S_{13}&&&+&2r_{45}\eta_2 && &=0,\\
	S_{14}&&&-&2r_{14}(\eta_2+\eta_3) & &        &=0,\\
	S_{15}&&&+&2r_{15}(\eta_2+\eta_3) & &        &=0,\\
	S_{23}&+ &2r_{45}\eta_1&&              &&&=0,\\
	S_{24}&-&2r_{24}(\eta_1-\eta_3)&&              & &        &=0,\\
	S_{25}&+&2r_{25}(\eta_1-\eta_3)& &              & &        &=0,\\
	S_{34}&+&2r_{34}(\eta_1+\eta_2) &&& &        &=0,\\
	S_{35}&-&2r_{35}(\eta_1+\eta_2) && & &        &=0,\\
	S_{45}&+&2r_{23}\eta_1&+&2r_{13}\eta_2 & +&2r_{12}\eta_3        &=0,\\
\end{matrix}
\right.
\end{equation}
where 
\begin{equation}\label{deltaSij}
\left\{
	\begin{aligned}
		S_{ij}=&m_im_jr_{ij}(\delta -r_{ij}^{-3}),\quad 1\leq i<j\leq n,\\
		\delta=&\frac{\lambda}{m}.
	\end{aligned}
\right.
\end{equation}
From \eqref{ten} we get 
\begin{equation}\label{etai}
\left\{
	\begin{aligned}
		\eta_1&=-\frac{S_{23}}{2r_{45}},\\
		\eta_2&=-\frac{S_{13}}{2r_{45}},\\
		\eta_3&=-\frac{S_{12}}{2r_{45}},\\
		\eta_2+\eta_3&=\frac{S_{14}}{2r_{14}}=-\frac{S_{15}}{2r_{15}},\\
		\eta_1-\eta_3&=\frac{S_{24}}{2r_{24}}=-\frac{S_{25}}{2r_{25}},\\
		\eta_1+\eta_2&=-\frac{S_{34}}{2r_{34}}=\frac{S_{35}}{2r_{35}},\\
		-r_{45}S_{45}&=r_{12}S_{12}+r_{13}S_{13}+r_{23}S_{23},
	\end{aligned}
	\right.
\end{equation}
where the last one is obtained by substituting the first three equations above to the last one in \eqref{ten}. We conclude directly from the last three equations above that 
\begin{equation}\label{i4i5}
	S_{i4}S_{i5}\leq0, \quad \forall i=1,2,3.
\end{equation}
The Hessian of $W_{123}(r)$ is
$D^2 W_{123}(r)=$
\tiny
\begin{equation*}
\begin{bmatrix}
R_{12} & 0 & 0 & 0 & 0 & 0 & 0 & 0 & 0 & 2\eta_3 \\
 0 & R_{13} & 0 & 0 & 0 & 0 & 0 & 0 & 0 & 2 \eta_2  \\
 0 & 0 & R_{14}\atop-2(\eta_2+\eta_3) & 0 & 0 & 0 & 0 & 0 & 0 & 0 \\
 0 & 0 & 0 & R_{15}\atop+2(\eta_2+\eta_3)  & 0 & 0 & 0 & 0 & 0 & 0 \\
 0 & 0 & 0 & 0 & R_{23} & 0 & 0 & 0 & 0 & 2\eta_1 \\
 0 & 0 & 0 & 0 & 0 & R_{24}\atop-2(\eta_1-\eta_3) & 0 & 0 & 0 & 0 \\
 0 & 0 & 0 & 0 & 0 & 0 & R_{25}\atop+2(\eta_1-\eta_3) & 0 & 0 & 0 \\
 0 & 0 & 0 & 0 & 0 & 0 & 0 & R_{34}\atop+2(\eta_1+\eta_2) & 0 & 0 \\
 0 & 0 & 0 & 0 & 0 & 0 & 0 & 0 & R_{35}\atop-2(\eta_1+\eta_2)  & 0 \\
 2\eta_3 & 2 \eta_2  & 0 & 0 & 2\eta_1 & 0 & 0 & 0 & 0 & R_{45} \\
\end{bmatrix},
\end{equation*}
\normalsize
where we denote by \begin{equation}\label{Rij}
	R_{ij}=m_im_j\left(\delta+\frac{2}{r_{ij}^3}\right)>0,\quad 1\leq i<j\leq 5.
\end{equation}
By direct computing, six positive eigenvalues of the above Hessian can be easily obtained with \eqref{etai}, and we denote by 
\begin{equation}
	\left\{
	\begin{aligned}
		\zeta_1&=R_{14}-2(\eta_2+\eta_3)=\frac{3}{r_{14}^3}>0,\\
		\zeta_2&=R_{15}+2(\eta_2+\eta_3)=\frac{3}{r_{15}^3}>0,\\
		\zeta_3&=R_{24}-2(\eta_1-\eta_3)=\frac{3}{r_{24}^3}>0,\\
		\zeta_4&=R_{25}+2(\eta_1-\eta_3)=\frac{3}{r_{25}^3}>0,\\
		\zeta_5&=R_{34}+2(\eta_1+\eta_2)=\frac{3}{r_{34}^3}>0,\\
		\zeta_6&=R_{35}-2(\eta_1+\eta_2)=\frac{3}{r_{35}^3}>0.\\
	\end{aligned}
	\right.
\end{equation}
The rest $\zeta_7,\cdots,\zeta_{10}$ are roots of the equation
\begin{equation}
	\mu^4+A_3\mu^3+A_2\mu^2+A_1\mu^1+A_0=0,
\end{equation}
where
\begin{equation}
	\left\{
	\begin{aligned}
		A_3=&-R_{12}-R_{13}-R_{23}-R_{45},\\
		A_2=&R_{12}R_{13}+R_{12}R_{23}+R_{13}R_{23}+(R_{12}+R_{13}+R_{23})R_{45}-4(\eta_1^2+\eta_2^2+\eta_3^2),\\
		A_1=&-R_{12}R_{13}R_{23}-(R_{12}R_{13}+R_{12}R_{23}+R_{13}R_{23})R_{45}\\
		      &+4\eta_1^2(R_{12}+R_{13})+4\eta_2^2(R_{12}+R_{23})+4\eta_3^2(R_{13}+R_{23}),\\
		A_0=&R_{12}R_{13}R_{23}R_{45}-4\eta_1^2R_{12}R_{13}-    4\eta_2^2R_{12}R_{23}- 4\eta_3^2R_{13}R_{23}.   
	\end{aligned}
	\right.
\end{equation}
If $\hat r$ is a critical point of $W_{123}(r)$, it satisfies \eqref{ten}. If we want to determine the type of this critical point, we need to check the sign of the ten eigenvalues of the Hessian at $\hat r$. Now we consider each sequential principle minors $f_i$ of the symmetric $D^2W_{123}(r)$ at $\hat r$ for $i=1,\cdots,10$. Since we can see obviously from the Hessian matrix that 
\begin{equation*}
	\left\{
	\begin{aligned}
	f_1=&R_{12}>0,\\
	f_2=&R_{12}R_{13}>0,\\
	f_3=&R_{12}R_{13}\zeta_1>0,\\
	f_4=&R_{12}R_{13}\zeta_1\zeta_2>0,\\
	f_5=&R_{12}R_{13}\zeta_1\zeta_2R_{23}>0,\\
	f_6=&R_{12}R_{13}\zeta_1\zeta_2R_{23}\zeta_3>0,\\
	f_7=&R_{12}R_{13}\zeta_1\zeta_2R_{23}\zeta_3\zeta_4>0,\\
	f_8=&R_{12}R_{13}\zeta_1\zeta_2R_{23}\zeta_3\zeta_4\zeta_5>0,\\
	f_9=&R_{12}R_{13}\zeta_1\zeta_2R_{23}\zeta_3\zeta_4\zeta_5>0,\\
	f_{10}=&(\zeta_1\cdots\zeta_6)\cdot \zeta_7\zeta_8\zeta_9\zeta_{10}.
	\end{aligned}
	\right.
\end{equation*}
Hence, it comes down to the sign of $f_{10}$, which comes down to the sign of $A_0=\zeta_7\zeta_8\zeta_9\zeta_{10}$.
To certify the positivity of $A_0$, we need to analyze the precise shape of this configuration with the equations in \eqref{etai} and the famous Perpendicular Bisector Theorem in \cite{moeckel1990} to get the relationships among $r_{ij}s$.
\subsection{Analyzing the precise shape} 
Without loss of generality, let $r_{12}\leq r_{23}$. For convenience, we denote by 
\begin{equation}\label{deltaij}
	\delta_{ij}=\delta-\frac{1}{r_{ij}^3}.
\end{equation}
Then the sign of $\delta_{ij}$ is consistent with that of $S_{ij}$. Furthermore we have $r_{ij}<r_{kl} \Leftrightarrow\delta_{ij}<\delta_{kl}$. In other words, compared with the positive $r_{ij}s$, $\delta_{ij}s$ are more subtle in showing the relationships whose signs indicate the difference with $\delta$.  

\textbf{By bisecting $\bf r_{12},r_{13},r_{15},r_{23},r_{34}$ and $\bf r_{45}$}, the following inequalities hold:
\begin{equation*}
\textcircled{0}
	\left\{
	\begin{aligned}
		\text{bisecting $r_{12}$: } r_{24}<r_{14} \text{ and } r_{15}<r_{25},\\
		\text{bisecting $r_{13}$: } r_{34}<r_{14} \text{ and } r_{15}<r_{35},\\
		\text{bisecting $r_{15}$: } r_{45}<r_{14} \text{ and } r_{12}<r_{25},\\
		\text{bisecting $r_{23}$: } r_{34}<r_{24} \text{ and } r_{25}<r_{35},\\
		\text{bisecting $r_{34}$: } r_{45}<r_{35} \text{ and } r_{23}<r_{24},\\
		\text{bisecting $r_{45}$: } r_{15}<r_{14} \text{ and } r_{34}<r_{35}.\\
	\end{aligned}
	\right.
\end{equation*}
We claim that $\delta_{12}\leq \delta_{23}<0$. Otherwise, suppose that $0\leq \delta_{12}\leq \delta_{23}$. Then from the third and the fifth lines in \textcircled{0} we have $\left\{\begin{aligned}
	&0\leq\delta_{12}<\delta_{25}\\
	&0\leq \delta_{23}<\delta_{24}
\end{aligned}\right.$, which contradicts with the fifth equation in \eqref{etai}. Furthermore we conclude from the sixth line in $\textcircled{0}$ and \eqref{i4i5} that
\begin{equation*}
	(\delta_{12},\delta_{23},)\delta_{15},\delta_{34}<0<\delta_{35},\delta_{14},
\end{equation*}
and the signs of the other four are derived from the following, in which the sign of $\delta_{13}$ and $\delta_{45}$ is of great importance.

\textbf{By bisecting $\bf r_{35}$}.
\begin{itemize}
	\item If \textcircled{1}$r_{23}<r_{25}$ holds, it is impossible to have $\left\{\begin{aligned}
		r_{13}\leq r_{15}\\
		r_{45}\leq r_{34}
	\end{aligned}\right.$, \\so \textcircled{3}$\left\{\begin{aligned}
		r_{13}\leq r_{15}\\
		r_{34}< r_{45}
	\end{aligned}\right.$ or \textcircled{4}$\left\{\begin{aligned}
		r_{15}< r_{13}\\
		r_{34}< r_{45}
	\end{aligned}\right.$ hold.
	\item If \textcircled{2}$r_{23}\geq r_{25}$ holds, it implies \textcircled{5} $r_{45}<r_{34}$.
\end{itemize}

\textbf{By bisecting $\bf r_{14}$}. It is impossible to have $\left\{\begin{aligned}
		r_{13}\leq r_{34}\\
		r_{45}\leq r_{15}
	\end{aligned}\right.$ since $r_{12}\leq r_{23}<r_{24}$. Hence \textcircled{6}$\left\{\begin{aligned}
		r_{13}\leq r_{34}\\
		r_{15}<r_{45}
	\end{aligned}\right.$ or \textcircled{7}$\left\{\begin{aligned}
		r_{34}<r_{13}\\
		r_{45}\leq r_{15}
	\end{aligned}\right.$ holds.
	
\textbf{By bisecting $\bf r_{24}$}. It is impossible to have $\left\{\begin{aligned}
		r_{25}\leq r_{45}\\
		r_{34}\leq r_{23}
	\end{aligned}\right.$ since $r_{12}<r_{14}$. Hence \textcircled{8}$\left\{\begin{aligned}
		r_{25}\leq r_{45}\\
		r_{23}< r_{34}
	\end{aligned}\right.$ or \textcircled{9}$\left\{\begin{aligned}
		r_{45}< r_{25}\\
		r_{34}\leq r_{23}
	\end{aligned}\right.$ holds.
 
\textbf{By bisecting $\bf r_{25}$}. It is impossible to have $\left\{\begin{aligned}
		r_{15}\leq  r_{12}\\
		r_{24}\leq r_{45}
	\end{aligned}\right.$ since $\delta_{23}<0<\delta_{35}$(, i.e., $r_{23}<r_{35}$). Hence \textcircled{10}$\left\{\begin{aligned}
		r_{12}< r_{15}\\
		r_{24}\leq r_{45}
	\end{aligned}\right.$ or \textcircled{11}$\left\{\begin{aligned}
		r_{15}\leq  r_{12}\\
		r_{45}< r_{24}
	\end{aligned}\right.$ holds.

\begin{lemma}\label{leqs}
With the given shape in Figure \ref{trapezoid} and the assumption of $r_{12}\leq r_{23}$, the mutual distances $r_{ij}s$ should satisfy the following inequalities
\begin{equation*}
	r_{12}\leq r_{23}<r_{13}\leq r_{34},r_{15}<r_{24},r_{25}\leq r_{45},r_{35},r_{14}.
\end{equation*}
More precisely, with the notation $\delta_{ij}$ defined in \eqref{deltaij} we have
\begin{equation}
\left\{\begin{aligned}
		\delta_{12}\leq \delta_{23}<\delta_{13}\leq \delta_{34}<\delta_{24}?&<0<\delta_{35},\delta_{14}\\
		\delta_{13}\leq \delta_{15}<\delta_{25}?&<0\leq\delta_{45}
	\end{aligned}\right.,
\end{equation}
where \enquote{?} denotes the signal is not determined.
\end{lemma}
\begin{proof}
From the \textcircled{0}\textcircled{1} path, we have a total of sixteen choices, and from \textcircled{0}\textcircled{2}, we have a total of eight. Meanwhile, some obvious contradictions can help us reduce the paths: 
\begin{itemize}
	\item \textcircled{3}-1 contradicts with \textcircled{11}-1, 
	\item \textcircled{6}-1 contradicts with \textcircled{9}-2, 
	\item \textcircled{7}-2 contradicts with \textcircled{8}-1, 
	\item \textcircled{4}-2 contradicts with \textcircled{10}-2;
	\item \textcircled{0}\textcircled{5} contradicts with \textcircled{10}-2.
\end{itemize}
This implies there are only four possible paths in \textcircled{0}\textcircled{1} and two in \textcircled{0}\textcircled{2} (where \enquote{?} denotes that the sign is uncertain): 
\begin{enumerate}[(i)]
	\item \label{i}\textcircled{0}\textcircled{1}\textcircled{3}\textcircled{6}\textcircled{8}\textcircled{10}: $\left\{\begin{aligned}
		\delta_{12}\leq \delta_{23}<\delta_{13}\leq \delta_{34}<\delta_{24}?&<0<\delta_{35},\delta_{14}\\
		\delta_{13}\leq \delta_{15}<\delta_{25}?&<0\leq\delta_{45}
	\end{aligned}\right.$;
	\item \label{ii}\textcircled{0}\textcircled{1}\textcircled{3}\textcircled{7}\textcircled{9}\textcircled{10}: $\left\{\begin{aligned}
		\delta_{34},\delta_{12}\leq \delta_{23}<\delta_{24}\leq \delta_{45}\leq \delta_{15}&<0<\delta_{25},\delta_{35},\delta_{14}\\
		 \delta_{23}<\delta_{13}\qquad\,\,\, \leq \delta_{15}&<0
	\end{aligned}\right.$;
	\item \label{iii}\textcircled{0}\textcircled{1}\textcircled{4}\textcircled{6}\textcircled{8}\textcircled{11}: $\left\{\begin{aligned}
		\delta_{15}\leq\delta_{12}\leq \delta_{23}<\delta_{25}\leq \delta_{45}\leq \delta_{34}&<0<\delta_{24},\delta_{35},\delta_{14}\\
		\delta_{23}<\delta_{13}\qquad\,\,\,\leq \delta_{34}&<0
	\end{aligned}\right.$;
	\item \label{iv}\textcircled{0}\textcircled{1}\textcircled{4}\textcircled{7}\textcircled{9}\textcircled{11}: $\left\{\begin{aligned}
		\delta_{45}\leq\delta_{15}\leq \delta_{12}\leq\delta_{23}&< \delta_{25}?,\delta_{24}?<0&<\delta_{35},\delta_{14}\\
		\delta_{45}\qquad\,\,\,\leq\delta_{34}\leq\delta_{23}&< \delta_{13}?\quad\quad\,<0&
	\end{aligned}\right.$;
	\item \label{v}\textcircled{0}\textcircled{2}\textcircled{5}\textcircled{6}\textcircled{8}\textcircled{10}: $\left\{\begin{aligned}
		\delta_{15}\leq \delta_{12}<\delta_{25}&\leq \delta_{23}<\delta_{13}\leq \delta_{34}&<0&<\delta_{24},\delta_{35},\delta_{14}\\
		\delta_{25}&\leq \delta_{45}\quad\quad\,\,\,<\delta_{34}&<0&
	\end{aligned}\right.$;
	\item \label{vi}\textcircled{0}\textcircled{2}\textcircled{5}\textcircled{7}\textcircled{9}\textcircled{11}: $\left\{\begin{aligned}
		\delta_{45}\leq\delta_{15}\leq \delta_{12}<\delta_{25}\leq \delta_{23}&<0<\delta_{24}&,\delta_{35},\delta_{14}\\
		\delta_{45}\,\qquad \qquad< \delta_{34}\quad\leq \delta_{23}&<\delta_{13}?,\delta_{24}&
	\end{aligned}\right.$.
\end{enumerate}
We claim that only (\ref{i}) is possible. 

For (\ref{ii}), $\delta_{34}<\delta_{15}$ implies that $P_5$ locates outside the segment $\overline{P_4'P_4''}$, where $P_1P_4'$ parallels $P_3P_4$ intersecting the line $P_4P_5$ at $P_4'$ and $P_4''$ is its axial symmetrical point with respect to the perpendicular line of $P_4P_5$ crossing $P_1$. Let $P_3P_5'$ be parallel to $P_1P_5$ intersecting $P_4P_5$ at $P_5'$. Without loss of generality, we assume that $P_4$ is located to the left of the foot of $P_3$'s perpendicular. 
\begin{itemize}
	\item See Figure \ref{0137910}. If $r_{13}>r_{45}$, then $P_5$ locates to the left of $P_4'$.  $P_5$ is located on the left side of $P_4$. By bisecting $r_{35}$, $r_{13}\leq r_{15}$ implies that $P_1$(resp. $P_5'$) locates to the left(resp. right) side of the intersection of the perpendicular bisector of $r_{35}$(the orange line) and $P_1P_3$(resp. $P_4P_5$). Since $r_{25}$ is located to the right side of $r_{35}$, the intersection of its perpendicular bisector(the pink line) and $P_4P_5$ is to the left side of that of $r_{35}$, i.e., it locates to the left side of $P_4$. This implies that $r_{24}>r_{45}$, from which we derive a contradiction.
	\item See Figure \ref{01379102}. Similarly, if $r_{13}\leq r_{45}$, then $P_5$ locates to the right side of $P_4''$. From $\delta_{34}\leq \delta_{13}$ we have the intersection of the perpendicular bisector of $r_{14}$(the green line) and $P_1P_3$(resp. $P_4P_5$) locates to the right(resp. left) side of $P_3$(resp. $P_4'$). This implies $r_{45}>r_{15}$, which contradicts with $\delta_{45}\leq \delta_{15}$.
\end{itemize}
\begin{figure}[htbp]
\centering
\subfigure[$r_{13}>r_{45}$.]{
\begin{minipage}{0.4\textwidth}
\centering
	\includegraphics[scale=0.33]{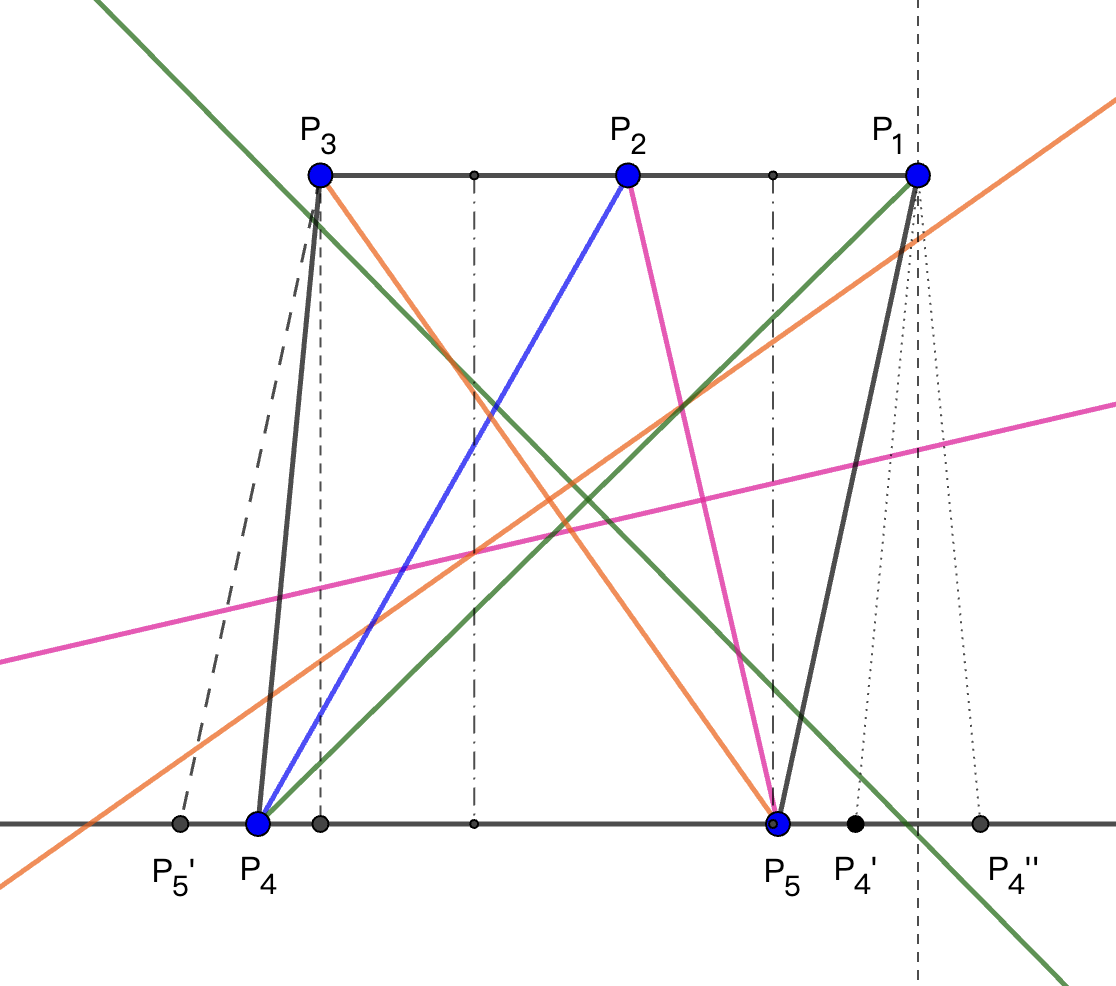}
\label{0137910}
\end{minipage}}
\qquad
\subfigure[$r_{13}\leq r_{45}$.]{
\begin{minipage}{0.4\textwidth}
\centering
\vspace{22pt}
	\includegraphics[scale=0.38]{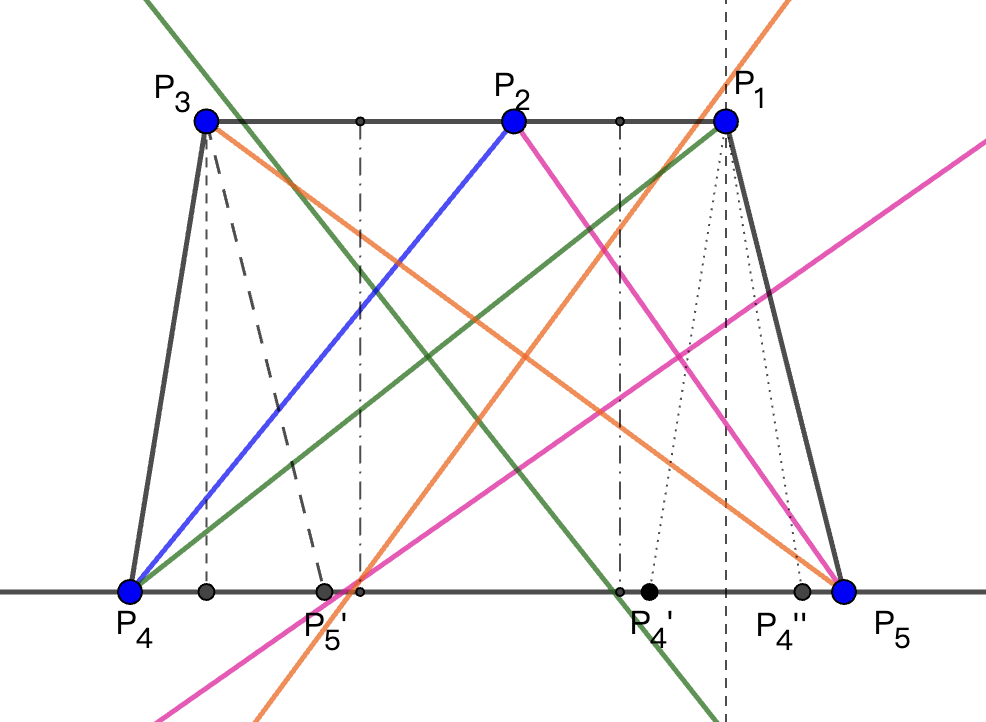}
\label{01379102}
\end{minipage}}
\caption{Contradictions in (\ref{ii}).}
\end{figure}

For (\ref{iii})-(\ref{vi}), the contradiction dues to that $r_{12}\geq r_{15}$ and $r_{24}>r_{45}$ cannot hold simultaneously. For example, for (\ref{iii}), see Figure \ref{0146811}. Let $P_2P_2'$ be parallel to $P_1P_5$ intersecting $P_4P_5$ at $P_2'$. Then $P_1P_2=P_5P_2'$. Let $Q_{23}, Q_{12}$ be the intersections of the perpendicular bisectors of $r_{23},r_{12}$ and $P_4P_5$ respectively. Then $P_4$ locates to the left side of $Q_{23}$ and $P_5$ to the right side of $Q_{12}$. Meanwhile, $P_2'$ locates to the right side of $Q_{23}$ since $P_5P_2'=P_1P_2\leq \frac{P_1P_3}{2}=Q_{12}Q_{23}$. From $r_{12}\geq r_{15}$ we conclude that $P_1$(resp. $P_2'$) locates to the right(resp. left) side of the intersection of $P_1P_2$(resp. $P_4P_5$) and the perpendicular bisector of $r_{25}$. Hence $P_4$ locates to the left side of the intersection on $P_4P_5$, which implies $r_{45}>r_{24}$, a contradiction. The other three cases are similar.
\begin{figure}
	\centering
	\includegraphics[scale=0.48]{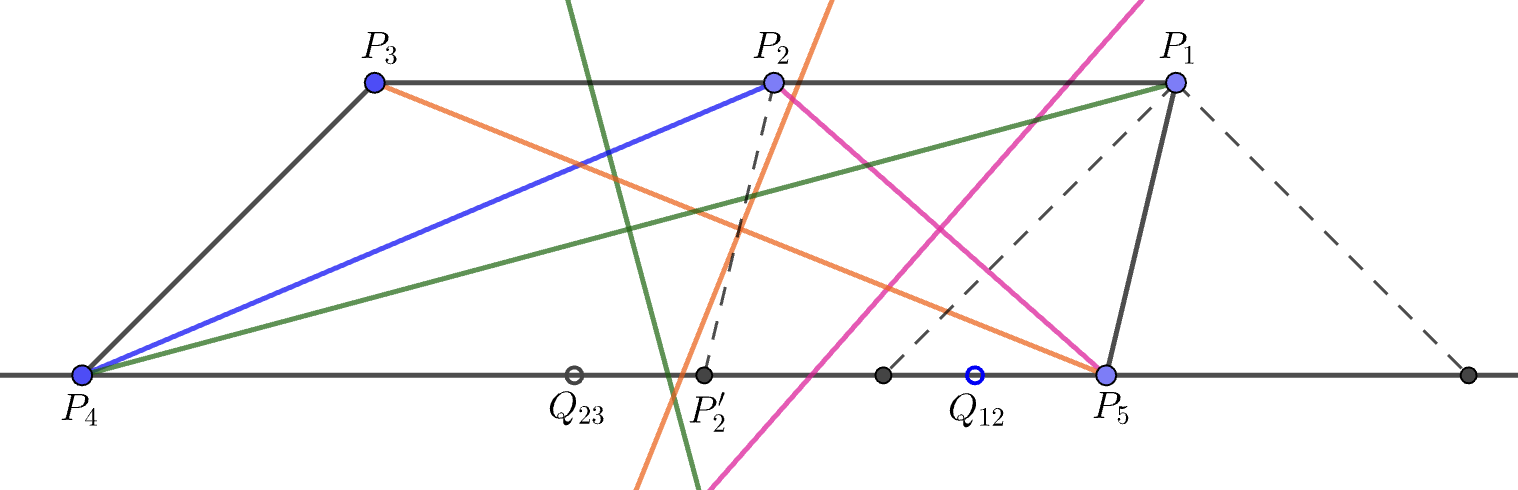}
	\caption{The contradiction in (\ref{iii}).}
	\label{0146811}
\end{figure}
\end{proof}

\subsection{The type of the critical point}
\begin{prop}\label{localmini}
Suppose that $\hat r\in \mathcal{M}^+_{123}$ is the critical point of $W_{123}(r)$. Then $\hat r$ is a non-degenerate local minimum of $W_{123}(r)$.
\end{prop}
\begin{proof}
	It is sufficient to show that $A_0>0$. 
	To simplify the calculation, we introduce
	\begin{equation}\label{Deltaij}
		\Delta_{ij}=\delta_{ij}+\frac{3}{r_{ij}^3}=\delta+\frac{2}{r_{ij}^3}>0.
	\end{equation}
	From the last equation in \eqref{etai} and the expressions in \eqref{deltaSij}, \eqref{deltaij} and \eqref{Rij} we have
	\begin{equation*}
	\left\{
	\begin{aligned}
		R_{45}&=-\frac{r_{12}S_{12}+r_{13}S_{13}+r_{23}S_{23}}{r_{45}^2}+\frac{3m_4m_5}{r_{45}^3},\\
		S_{ij}&=r_{ij}m_im_j\left(\Delta_{ij}-\frac{3}{r_{ij}^3}\right),\\
		R_{ij}&=m_im_j\Delta_{ij},\quad \forall i\neq j.	
	\end{aligned}
		\right.
	\end{equation*}
Substituting the above formulas into $A_0$ sequentially(i.e., first $R_{45}$, then $R_{12},R_{13}$ and $R_{23}$) and we will obtain an expression in terms of $\Delta_{12},\Delta_{13}$ and $\Delta_{23}$:
\begin{equation}\label{A0}
\begin{aligned}
	A_0=-\frac{m_1^2m_2^2m_3^2}{r_{12}^4r_{13}^4r_{23}^4r_{45}^3}(&-3m_4m_5r_{12}^4r_{13}^4r_{23}^4\Delta_{12}\Delta_{13}\Delta_{23}\\
	&+r_{45}m_1m_3r_{12}^4r_{23}^4\Delta_{12}\Delta_{23}(9-9r_{13}^3\Delta_{13}+2r_{13}^6\Delta_{13}^2)\\
	&+r_{45}m_1m_2r_{13}^4r_{23}^4\Delta_{13}\Delta_{23}(9-9r_{12}^3\Delta_{12}+2r_{12}^6\Delta_{12}^2)\\
	&+r_{45}m_2m_3r_{12}^4r_{13}^4\Delta_{12}\Delta_{13}(9-9r_{23}^3\Delta_{23}+2r_{23}^6\Delta_{23}^2)).
\end{aligned}
\end{equation}
Then for $i\neq j$ and $i,j \in\left\{1,2,3\right\}$ we have
\begin{equation*}
\begin{aligned}
9-9r_{ij}^3\Delta_{ij}+2r_{ij}^6\Delta_{ij}^2
&=2r_{ij}^6\left(\Delta_{ij}-\frac{9}{4r_{ij}^3}\right)^2-\frac{9}{8}\\
&=2r_{ij}^6\left(\delta_{ij}-\frac{3}{4r_{ij}^3}\right)^2-\frac{9}{8}<0.
\end{aligned}
\end{equation*}
We notice that the corresponding $\delta_{12},\delta_{13},\delta_{23}$ are all negative from  Lemma \ref{leqs}, which helps us to determine $A_0>0$.
\end{proof}

\subsection{The topology of $\mathcal{M}_{123}^+(=\mathcal{M}_{T_1T_3L}^+)$}
Let $$\mathcal{M}_{T_1T_3L}=\left\{r\in(\mathbb{R}^+)^{10}\vert I(r)=I_0,T_1(r)=0, T_2(r)=0 ,T_{3}(r)=0\right\},$$
and  let $\mathcal{M}_{T_1T_3L}^0$ and $\mathcal{M}_{T_1T_3L}^{0+}$ be the sets corresponding to the equal mass case of $\mathcal{M}_{T_1T_3L}$ and $\mathcal{M}_{T_1T_3L}^+$ respectively. 
\begin{lemma}
	$\chi(\mathcal{M}_{123}^+)=1$.
\end{lemma}
\begin{proof}
Without loss of generality, let $m_1=\cdots=m_5=1$. Then for $r\in (\mathbb{R}^+)^{10}$ the constraints in $\mathcal{M}_{T_1T_3L}^0$ are 
\begin{equation}
\begin{aligned}
 r_{12}^2+r_{13}^2+r_{14}^2+r_{15}^2+r_{23}^2+r_{24}^2+r_{25}^2+r_{34}^2+r_{35}^2+r_{45}^2&=1,\\
 T_1(r)=2r_{23}r_{45}-r_{24}^2+r_{25}^2+r_{34}^2-r_{35}^2&=0,\\
T_3(r)=2r_{12}r_{45}-r_{14}^2+r_{15}^2+r_{24}^2-r_{25}^2&=0,\\
L_{1,2,3}=r_{12}-r_{13}+r_{23}&=0.
\end{aligned}
\end{equation}
From the first and the second equations, we have 
\begin{equation}\label{t1}
\left\{
	\begin{aligned}
	(r_{23}+r_{45})^2+r_{12}^2+r_{14}^2+r_{15}^2+0r_{23}^2+0r_{24}^2+2r_{25}^2+2r_{34}^2+0r_{35}^2&=1-(r_{13}^2)	,\\
	(r_{23}-r_{45})^2+r_{12}^2+r_{14}^2+r_{15}^2+0r_{23}^2+2r_{24}^2+0r_{25}^2+0r_{34}^2+2r_{35}^2&=1-(r_{13}^2)	,\\
	\end{aligned}
	\right.
\end{equation}
and similarly, from the first and the third equations, we have 
\begin{equation}\label{t2}
\left\{
	\begin{aligned}
	(r_{12}+r_{45})^2+0r_{12}^2+0r_{14}^2+2r_{15}^2+r_{23}^2+2r_{24}^2+0r_{25}^2+r_{34}^2+r_{35}^2&=1-(r_{13}^2)	,\\
	(r_{12}-r_{45})^2+0r_{12}^2+2r_{14}^2+0r_{15}^2+r_{23}^2+0r_{24}^2+2r_{25}^2+r_{34}^2+r_{35}^2&=1-(r_{13}^2).\\		\end{aligned}
	\right.
\end{equation}
Now let 
\begin{equation*}
\left\{
\begin{aligned}
r_{23}+r_{45}&=v_1,\\
r_{23}-r_{45}&=w_1,
\end{aligned}
\right.
\end{equation*}
which implies 
\begin{equation}\label{v1w1}
\left\{
\begin{aligned}
r_{23}&=\frac{v_1+w_1}{2}\geq0,\\
r_{45}&=\frac{v_1-w_1}{2}\geq0,\\
r_{12}&=r_{13}-\frac{v_1+w_1}{2}\geq0.
\end{aligned}
\right.
\end{equation}
Let $r_{13}^2+r_{12}^2+r_{14}^2+r_{15}^2=\varepsilon\in [0,1]$. Then \eqref{t1} becomes 
\begin{equation}
\left\{
	\begin{aligned}
		v_1^2+2r_{25}^2+2r_{34}^2&=1-\varepsilon,\\
		w_1^2+2r_{24}^2+2r_{35}^2&=1-\varepsilon.\\
	\end{aligned}
	\right.
\end{equation}
Denote by 
\begin{equation*}
	\begin{aligned}
		\mathcal{E}_1^2(1-\varepsilon)&=\left\{(v_1,r_{25},r_{34})\in \mathbb{R}^3\vert  v_1^2+2r_{25}^2+2r_{34}^2=1-\varepsilon\right\},\\
		\mathcal{E}_2^2(1-\varepsilon)&=\left\{(w_1,r_{24},r_{35})\in \mathbb{R}^3\vert w_1^2+2r_{24}^2+2r_{35}^2=1-\varepsilon\right\}.
	\end{aligned}
\end{equation*}

Firstly, with $\varepsilon$ fixed, from \eqref{t1}, the projection 
$$\tilde p_1: (v_1,r_{25},r_{34},w_1,r_{24},r_{35})\to (v_1,r_{25},r_{34})$$ 
induces a contractible fibration $p_1$, since the base $$\left\{(v_1,r_{25},r_{34})\in \mathcal{E}_1^2(1-\varepsilon)\vert v_1,r_{25},r_{34}\geq 0\right\}$$ and the fiber $$\left\{(w_1,r_{24},r_{35})\in \mathcal{E}_2^2(1-\varepsilon)\vert -\sqrt{1-\varepsilon}\leq w_1\leq v_1,r_{24},r_{35}\geq0 \right\}$$ are both contractible. 

Secondly, we consider another projection 
$$\tilde p_2:(v_1,r_{25},r_{34},w_1,r_{24},r_{35},r_{12},r_{13},r_{14},r_{15})\to (v_1,r_{25},r_{34},w_1,r_{24},r_{35}).$$ 
For any fixed $(v_1,r_{25},r_{34},w_1,r_{24},r_{35})$ in the contractible fibration in the above $p_1$, the region formed by $(r_{12},r_{13},r_{14},r_{15})$ is contractible. 
In fact, noticing that if $r_{13}$ is fixed in $[\frac{v_1+w_1}{2},\varepsilon]$, then $r_{12}$ and $r_{23}$ are uniquely determined with the relationship in \eqref{v1w1}. Hence, $r_{14}$ and $r_{15}$ are uniquely determined by the second and the first equation in \eqref{t2}, respectively. 
Hence the region formed by $(r_{12},r_{13},r_{14},r_{15})$ is homeomorphic to a contractible segment of $r_{13}$. 
Therefore, a new contractible fibration $p_2$ can be induced from $\tilde p_2$ with both the base and the fiber contractible. This gives the contractibility of $\mathcal{M}^{0+}_{123}$. 

Since $I(r)=I_0$ forms a nine-dimensional ellipsoid $\mathcal{E}^9$, which is homeomorphic to $\mathcal{S}^9$, i.e., the nine-dimensional sphere with $m_1=\cdots=m_5=1$. Meanwhile the cones $T_1,T_3=0$ and the hyperplane $L_{1,2,3}=0$ passing through the origin intersect $\mathcal{E}^9_0$ in $\mathcal{M}^{0+}_{T_1T_3L}\subset \mathcal{E}^9_0$ corresponding to $\mathcal{M}^{+}_{T_1T_3L}$ via the homeomorphism. Hence we have $\chi(\mathcal{M}^{0+}_{123})=\chi(\mathcal{M}^{0+}_{T_1T_3L})=\chi(\mathcal{M}^{+}_{T_1T_3L})=1$.
\end{proof}
And so 
\begin{lemma}
	The potential $U$ has a unique critical point on $\mathcal{M}_{123}^+$.
\end{lemma}
\begin{proof}
From Lemma \ref{localmini}, critical points of $U$ restricting on $\mathcal{M}_{123}^+$ are non-degenerate local minimum with Morse index 0, which implies that $U$ is a Morse function. Noticing that when $r\to \partial\mathcal{M}_{123}^+$ we have $U\to+\infty$, so the minimum exists. Then, from the Morse inequality 
$$\sum (-1)^q\alpha_q=\chi (\mathcal{M}_{123}^{+})=1$$
we have $\alpha_0=1$ provided $q=0$. 
\end{proof}

Now, we have reached the final step:
\begin{proof}[The proof of Theorem \ref{5planarunique}]
From Proposition \ref{rinH123iffx}, the distance vector $r$ of a central configuration satisfying the above condition is the critical point of the Lagrangian function $W_{123}$, namely, the potential function $U$ restricted on $\mathcal{M}_{123}^+$. From Lemma \ref{localmini}, this critical point must be a non-degenerate local minimum. Since $\mathcal{T}_{123}\subset \mathcal{M}_{123}^+$, we conclude that $U$ has at most one critical point restricted on $\mathcal{T}_{123}$.
\end{proof}

\subsection{The proof of Proposition \texorpdfstring{\ref{5co}}{Lg}}\label{proof5co}

By direct computation, we have $F_i(x)\equiv0$; hence, the first "if and only if" holds. Secondly, noticing that $r=r(x)$, by the chain rule
\begin{equation*}
		\nabla_x \tilde W(r(x))=\nabla_r \tilde W(r)\cdot\nabla_x r(x), 
	\end{equation*}
it is easy to see that if $r$ is a critical point of $\tilde W(r)$, then $x$ is a critical point of $\tilde W(r(x))$. Now suppose that $x$ is a critical point of $\tilde W(r(x))$. With the notations in \eqref{deltaSij} and direct computation, we write down the components of $\nabla_x \tilde W(r(x))$
\begin{subequations}
\begin{align}
  \frac{\partial \tilde W(r(x))}{\partial x_1}&=\frac{S_{12}}{r_{12}}\overrightarrow{P_2P_1}+\frac{S_{13}}{r_{13}}\overrightarrow{P_3P_1}+\frac{S_{14}}{r_{14}}\overrightarrow{P_4P_1}+\frac{S_{15}}{r_{15}}\overrightarrow{P_5P_1}=0,\label{x1}\\
  \frac{\partial \tilde W(r(x))}{\partial x_2}&=\frac{S_{12}}{r_{12}}\overrightarrow{P_1P_2}+\frac{S_{23}}{r_{23}}\overrightarrow{P_3P_2}+\frac{S_{24}}{r_{24}}\overrightarrow{P_4P_2}+\frac{S_{25}}{r_{25}}\overrightarrow{P_5P_2}=0,\label{x2}\\
  \frac{\partial \tilde W(r(x))}{\partial x_3}&=\frac{S_{13}}{r_{13}}\overrightarrow{P_1P_3}+\frac{S_{23}}{r_{23}}\overrightarrow{P_2P_3}+\frac{S_{34}}{r_{34}}\overrightarrow{P_4P_3}+\frac{S_{35}}{r_{35}}\overrightarrow{P_5P_3}=0,\label{x3}\\
\frac{\partial \tilde W(r(x))}{\partial x_4}&=\frac{S_{14}}{r_{14}}\overrightarrow{P_1P_4}+\frac{S_{24}}{r_{24}}\overrightarrow{P_2P_4}+\frac{S_{34}}{r_{34}}\overrightarrow{P_3P_4}+\frac{S_{45}}{r_{45}}\overrightarrow{P_5P_4}=0,\label{x4}\\
\frac{\partial \tilde W(r(x))}{\partial x_5}&=\frac{S_{15}}{r_{15}}\overrightarrow{P_1P_5}+\frac{S_{25}}{r_{25}}\overrightarrow{P_2P_5}+\frac{S_{35}}{r_{35}}\overrightarrow{P_3P_5}+\frac{S_{45}}{r_{45}}\overrightarrow{P_4P_5}=0,\label{x5}
\end{align}
\end{subequations}
where $\overrightarrow{P_iP_j}$ is the vector from $P_i$ to $P_j$. Meanwhile noticing that for the Cayley-Menger determinant $\mathcal{F}_4(i,j,k,l)$ formed by four particles $P_i,P_j,P_k$ and $P_l$ we have 
\begin{equation}\label{F4partialrij}
	\frac{\partial \mathcal{F}_4(i,j,k,l)}{\partial r_{ij}^2}=-32 \triangle_i\triangle_j,
\end{equation}
where $\triangle_i=\triangle_{j,k,l}$ represents the oriented area of the triangle formed by the other particles with the $i$th deleted, whose sign is determined by the position $t$ of $i$ in “$i,j,k,l$,” i.e., $(-1)^{t+1}$. In other words, the value of $\triangle_i$ is the “algebra minor” 
\begin{equation*}
	\triangle_{i}=\triangle_{j,k,l}=\frac{(-1)^{t+1}}{2}\begin{vmatrix}
		1&1&1\\
		x_j &x_k&x_l
	\end{vmatrix}=\frac{(-1)^{t+1}}{2}\overrightarrow{P_lP_j}\times\overrightarrow{P_lP_k},
\end{equation*}
by deleting the $t$th column of the $i$th body in the matrix 
$$\begin{bmatrix}
		1&1&1&1\\
		x_i&x_j &x_k&x_l
	\end{bmatrix},$$ 
where “$\times$” represents the wedge product of two vectors. It should be pointed out that with the order $i<j<k<l$, if the sign of the oriented area is "$+$", it coincides with the wedge product of the corresponding vectors; if the sign is "$-$", it differs from the wedge product only by a negative sign. In the following, we use $\triangle_i^+$ and $\triangle_j^-$ to represent the signs of the oriented areas for convenience directly. 

From the above, we have 
\begin{equation}\label{F4/rij}
\frac{\partial \mathcal{F}_4(i,j,k,l)}{\partial r_{ij}}=	\frac{\partial \mathcal{F}_4(i,j,k,l)}{\partial r_{ij}^2}\cdot \frac{\partial r_{ij}^2}{\partial r_{ij}}=2r_{ij}\cdot(-32 \triangle_i\triangle_j)=-64r_{ij}\triangle_i\triangle_j.
\end{equation}
Now we write down the components of $\nabla_r \tilde W(r)$ 
\begin{subequations}
\begin{align}
	\frac{\partial \tilde W}{\partial r_{12}}=&S_{12}-64r_{12}(\eta_5\triangle_{234}^{+}\triangle_{134}^{-}+\eta_4\triangle_{235}^{+}\triangle_{135}^{-}+\eta_3\triangle_{245}^{+}\triangle_{145}^{-}),\label{12}\\
	\frac{\partial \tilde W}{\partial r_{13}}=&S_{13}-64r_{13}(\eta_5\triangle_{234}^{+}\triangle_{124}^{+}+\eta_4\triangle_{235}^{+}\triangle_{125}^{+}),\label{13}\\
	\frac{\partial \tilde W}{\partial r_{14}}=&S_{14}-64r_{14}(\eta_5\triangle_{234}^{+}\triangle_{123}^{-}+\eta_3\triangle_{245}^{+}\triangle_{125}^{+}),\label{14}\\
	\frac{\partial \tilde W}{\partial r_{15}}=&S_{15}-64r_{15}(\eta_4\triangle_{235}^{+}\triangle_{123}^{-}+\eta_3\triangle_{245}^{+}\triangle_{124}^{-}),\label{15}\\
	\frac{\partial \tilde W}{\partial r_{23}}=&S_{23}-64r_{23}(\eta_5\triangle_{134}^{-}\triangle_{124}^{+}+\eta_4\triangle_{135}^{-}\triangle_{125}^{+}),\label{23}\\
	\frac{\partial \tilde W}{\partial r_{24}}=&S_{24}-64r_{24}(\eta_5\triangle_{134}^{-}\triangle_{123}^{-}+\eta_3\triangle_{145}^{-}\triangle_{125}^{+}),\label{24}\\
	\frac{\partial \tilde W}{\partial r_{25}}=&S_{25}-64r_{25}(\eta_4\triangle_{135}^{-}\triangle_{123}^{-}+\eta_3\triangle_{145}^{-}\triangle_{124}^{-}),\label{25}\\
	\frac{\partial \tilde W}{\partial r_{34}}=&S_{34}-64r_{34}\eta_5\triangle_{124}^{+}\triangle_{123}^{-},\label{34}\\
	\frac{\partial \tilde W}{\partial r_{35}}=&S_{35}-64r_{35}\eta_4\triangle_{125}^{+}\triangle_{123}^{-},\label{35}\\
	\frac{\partial \tilde W}{\partial r_{45}}=&S_{45}-64r_{45}\eta_3\triangle_{125}^{+}\triangle_{124}^{-}.\label{45}
\end{align}
\end{subequations}
The conditions on $F_i$ guarantee all the $\triangle$s appearing in the denominators not to be zero.
We conclude that when $\nabla_x\tilde W(r(x))=0$, i.e., \eqref{x1}-\eqref{x5} hold, then \eqref{12}-\eqref{25} hold provided 
\begin{equation}\label{eta345}
\left\{
\begin{aligned}
	\eta_5=\frac{1}{64}\frac{S_{34}}{r_{34}}\frac{1}{\triangle_{124}^{+}\triangle_{123}^{-}},\\
	\eta_4=\frac{1}{64}\frac{S_{35}}{r_{35}}\frac{1}{\triangle_{125}^{+}\triangle_{123}^{-}},\\
	\eta_3=\frac{1}{64}\frac{S_{45}}{r_{45}}\frac{1}{\triangle_{125}^{+}\triangle_{124}^{-}},\\
\end{aligned}
\right.
\end{equation}
namely 
\eqref{34},\eqref{35} and \eqref{45} hold.

\begin{enumerate}[(i)]
	\item For \eqref{25}, substituting \eqref{eta345} into the following
\begin{equation}\label{W/r25}
\begin{aligned}
	\frac{\partial \tilde W}{\partial r_{25}}=&S_{25}-64r_{25}(\eta_4\triangle_{135}^{-}\triangle_{123}^{-}+\eta_3\triangle_{145}^{-}\triangle_{124}^{-})\\
   =&\frac{r_{25}}{\triangle_{125}^{+}}\left(
   \frac{S_{25}}{r_{25}}\triangle_{125}^{+}
   -\frac{S_{35}}{r_{35}}\triangle_{135}^{-}
   -\frac{S_{45}}{r_{45}}\triangle_{145}^{-}\right)\\
   =&\frac{r_{25}}{\triangle_{125}}\left(\frac{S_{25}}{r_{25}}\overrightarrow{P_5P_1}\times \overrightarrow{P_5P_2} -\frac{S_{35}}{r_{35}}\overrightarrow{P_5P_1}\times \overrightarrow{P_5P_3}-\frac{S_{45}}{r_{45}}\overrightarrow{P_5P_1}\times \overrightarrow{P_5P_4}\right)\\
   =&\frac{r_{25}}{\triangle_{125}}\overrightarrow{P_5P_1}\times  \left(\frac{S_{25}}{r_{25}} \overrightarrow{P_5P_2} 
   +\frac{S_{35}}{r_{35}} \overrightarrow{P_5P_3}
   +\frac{S_{45}}{r_{45}}\overrightarrow{P_5P_4}\right),
\end{aligned}
\end{equation}
and from \eqref{x5}, the formula in the parenthesis of the last line of \eqref{W/r25} can be used to compute
$$\frac{S_{25}}{r_{25}} \overrightarrow{P_5P_2} 
   +\frac{S_{35}}{r_{35}} \overrightarrow{P_5P_3}
   +\frac{S_{45}}{r_{45}}\overrightarrow{P_5P_4}=\frac{S_{15}}{r_{15}} \overrightarrow{P_1P_5}.$$
Noticing that $\overrightarrow{P_5P_1}\times\overrightarrow{P_1P_5}=0$, hence $\frac{\partial \tilde W}{\partial r_{25}}=0$. Similarly, 

\item From \eqref{24} and \eqref{x4}, noticing that
$$\frac{S_{24}}{r_{24}}\overrightarrow{P_4P_2}
+\frac{S_{34}}{r_{34}}\overrightarrow{P_4P_3}
+\frac{S_{45}}{r_{45}}\overrightarrow{P_4P_5}
=\frac{S_{14}}{r_{14}}\overrightarrow{P_1P_4},$$
we have 
\begin{equation}\label{W/r24}
\begin{aligned}
\frac{\partial \tilde W}{\partial r_{24}}=&S_{24}-64r_{24}(\eta_5\triangle_{134}^{-}\triangle_{123}^{-}+\eta_3\triangle_{145}^{-}\triangle_{125}^{+})\\
=&S_{24}-r_{24}\left(\frac{S_{34}}{r_{34}}\frac{\triangle_{134}^{-}}{\triangle_{124}^{+}}
+\frac{S_{45}}{r_{45}}\frac{\triangle_{145}^{-}}{\triangle_{124}^{-}}\right)\\
=&\frac{r_{24}}{\triangle_{124}^{+}}\left(
\frac{S_{24}}{r_{24}}\triangle_{124}^{+}
-\frac{S_{34}}{r_{34}}\triangle_{134}^{-}
-\frac{S_{45}}{r_{45}}\triangle_{145}^{-}
\right)\\
=&\frac{r_{24}}{\triangle_{124}^{+}}\left(
\frac{S_{24}}{r_{24}}\overrightarrow{P_4P_1}\times \overrightarrow{P_4P_2}
-\frac{S_{34}}{r_{34}}(-\overrightarrow{P_4P_1}\times \overrightarrow{P_4P_3})
-\frac{S_{45}}{r_{45}}(-\overrightarrow{P_4P_1}\times \overrightarrow{P_4P_5})\right)\\
=&\frac{r_{24}}{\triangle_{124}^{+}}\overrightarrow{P_4P_1}\times
\left(
\frac{S_{24}}{r_{24}}\overrightarrow{P_4P_2}
+\frac{S_{34}}{r_{34}}\overrightarrow{P_4P_3}
+\frac{S_{45}}{r_{45}}\overrightarrow{P_4P_5}
\right)\\
=&0.
\end{aligned}
\end{equation}

\item From \eqref{23} and \eqref{x3}, noticing that
$$\frac{S_{23}}{r_{23}}\overrightarrow{P_3P_2}
+\frac{S_{34}}{r_{34}}\overrightarrow{P_3P_4}
+\frac{S_{35}}{r_{35}}\overrightarrow{P_3P_5}=\frac{S_{13}}{r_{13}}\overrightarrow{P_1P_3},$$
we have 
\begin{equation*}
\begin{aligned}
\frac{\partial \tilde W}{\partial r_{23}}=&S_{23}-64r_{23}(\eta_5\triangle_{134}^{-}\triangle_{124}^{+}+\eta_4\triangle_{135}^{-}\triangle_{125}^{+})\\
=&S_{23}-r_{23}\left(\frac{S_{34}}{r_{34}}\frac{\triangle_{134}^{-}}{\triangle_{123}^{-}}
+\frac{S_{35}}{r_{35}}\frac{\triangle_{135}^{-}}{\triangle_{123}^{-}}\right)\\
=&\frac{r_{23}}{\triangle_{123}^{-}}
\left(
\frac{S_{23}}{r_{23}}\triangle_{123}^{-}
-\frac{S_{34}}{r_{34}}\triangle_{134}^{-}
-\frac{S_{35}}{r_{35}}\triangle_{135}^{-}
\right)\\
=&\frac{r_{23}}{\triangle_{123}^{-}}
\left(
\frac{S_{23}}{r_{23}}\triangle_{123}^{-}
+\frac{S_{34}}{r_{34}}\triangle_{143}^{-}
+\frac{S_{35}}{r_{35}}\triangle_{153}^{-}
\right)\\
=&\frac{r_{23}}{\triangle_{123}^{-}}\overrightarrow{P_3P_1}\times 
\left(
-\frac{S_{23}}{r_{23}}\overrightarrow{P_3P_2}
-\frac{S_{34}}{r_{34}}\overrightarrow{P_3P_4}
-\frac{S_{35}}{r_{35}}\overrightarrow{P_3P_5}
\right)\\
=&0.
\end{aligned}
\end{equation*}

\item From \eqref{15} and \eqref{x5}, noticing that
$$\frac{S_{15}}{r_{15}}\overrightarrow{P_5P_1}
+\frac{S_{35}}{r_{35}}\overrightarrow{P_5P_3}
+\frac{S_{45}}{r_{45}}\overrightarrow{P_5P_4}
=\frac{S_{25}}{r_{25}}\overrightarrow{P_2P_5},
$$
we have
\begin{equation}\label{W/r15}
	\begin{aligned}
\frac{\partial \tilde W}{\partial r_{15}}
=&S_{15}-64r_{15}(\eta_4\triangle_{235}^{+}\triangle_{123}^{-}+\eta_3\triangle_{245}^{+}\triangle_{124}^{-})\\
=&S_{15}-r_{15}\left(
\frac{S_{35}}{r_{35}}\frac{\triangle_{235}^{+}}{\triangle_{125}^{+}}
+\frac{S_{45}}{r_{45}}\frac{\triangle_{245}^{+}}{\triangle_{125}^{+}}\right)\\
=&\frac{r_{15}}{\triangle_{125}^{+}}
\left(
\frac{S_{15}}{r_{15}}\triangle_{125}^{+}
-\frac{S_{35}}{r_{35}}\triangle_{235}^{+}
-\frac{S_{45}}{r_{45}}\triangle_{245}^{+}
\right)\\
=&\frac{r_{15}}{\triangle_{125}^{+}}
\left(
-\frac{S_{15}}{r_{15}}\triangle_{215}^{+}
-\frac{S_{35}}{r_{35}}\triangle_{235}^{+}
-\frac{S_{45}}{r_{45}}\triangle_{245}^{+}
\right)\\
=&\frac{r_{15}}{\triangle_{125}^{+}}\overrightarrow{P_5P_2}\times 
\left(
-\frac{S_{15}}{r_{15}}\overrightarrow{P_5P_1}
-\frac{S_{35}}{r_{35}}\overrightarrow{P_5P_3}
-\frac{S_{45}}{r_{45}}\overrightarrow{P_5P_4}
\right)\\
=&0.
\end{aligned}
\end{equation}

\item From \eqref{14} and \eqref{x4}, noticing that
$$\frac{S_{14}}{r_{14}}\overrightarrow{P_4P_1}
+\frac{S_{34}}{r_{34}}\overrightarrow{P_4P_3}
+\frac{S_{45}}{r_{45}}\overrightarrow{P_4P_5}
=\frac{S_{24}}{r_{24}}\overrightarrow{P_2P_4},
$$
we have 
\begin{equation}\label{W/r14}
\begin{aligned}
\frac{\partial \tilde W}{\partial r_{14}}
=&S_{14}-64r_{14}(\eta_5\triangle_{234}^{+}\triangle_{123}^{-}+\eta_3\triangle_{245}^{+}\triangle_{125}^{+})\\
=&S_{14}-r_{14}\left(
\frac{S_{34}}{r_{34}}\frac{\triangle_{234}^{+}}{\triangle_{124}^{+}}
+\frac{S_{45}}{r_{45}}\frac{\triangle_{245}^{+}}{\triangle_{124}^{-}}
\right)\\
=&\frac{r_{14}}{\triangle_{124}^{+}}
\left(
\frac{S_{14}}{r_{14}}\triangle_{124}^{+}
-\frac{S_{34}}{r_{34}}\triangle_{234}^{+}
+\frac{S_{45}}{r_{45}}\triangle_{245}^{+}
\right)\\
=&\frac{r_{14}}{\triangle_{124}^{+}}
\left(
-\frac{S_{14}}{r_{14}}\triangle_{214}^{+}
-\frac{S_{34}}{r_{34}}\triangle_{234}^{+}
-\frac{S_{45}}{r_{45}}\triangle_{254}^{+}
\right)\\
=&\frac{r_{14}}{\triangle_{124}^{+}}\overrightarrow{P_4P_2}\times 
\left(
-\frac{S_{14}}{r_{14}}\overrightarrow{P_4P_1}
-\frac{S_{34}}{r_{34}}\overrightarrow{P_4P_3}
-\frac{S_{45}}{r_{45}}\overrightarrow{P_4P_5}
\right)\\
=&0.
\end{aligned}
\end{equation}

\item From \eqref{13} and \eqref{x3}, noticing that
$$\frac{S_{13}}{r_{13}}\overrightarrow{P_3P_1}
+\frac{S_{34}}{r_{34}}\overrightarrow{P_3P_4}
+\frac{S_{35}}{r_{35}}\overrightarrow{P_3P_5}
=\frac{S_{23}}{r_{23}}\overrightarrow{P_2P_3},
$$
we have
\begin{equation*}
\begin{aligned}
\frac{\partial \tilde W}{\partial r_{13}}
=&S_{13}-64r_{13}(\eta_5\triangle_{234}^{+}\triangle_{124}^{+}+\eta_4\triangle_{235}^{+}\triangle_{125}^{+})\\
=&S_{13}-r_{13}\left(
\frac{S_{34}}{r_{34}}\frac{\triangle_{234}^{+}}{\triangle_{123}^{-}}
+\frac{S_{35}}{r_{35}}\frac{\triangle_{235}^{+}}{\triangle_{123}^{-}}
\right)\\
=&\frac{r_{13}}{\triangle_{123}^{-}}
\left(
\frac{S_{13}}{r_{13}}\triangle_{123}^{-}
-\frac{S_{34}}{r_{34}}\triangle_{234}^{+}
-\frac{S_{35}}{r_{35}}\triangle_{235}^{+}
\right)\\
=&\frac{r_{13}}{\triangle_{123}^{-}}
\left(
-\frac{S_{13}}{r_{13}}\triangle_{213}^{-}
+\frac{S_{34}}{r_{34}}\triangle_{243}^{+}
+\frac{S_{35}}{r_{35}}\triangle_{253}^{+}
\right)\\
=&\frac{r_{13}}{\triangle_{123}^{-}}\overrightarrow{P_3P_2}\times 
\left(
\frac{S_{13}}{r_{13}}\overrightarrow{P_3P_1}
+\frac{S_{34}}{r_{34}}\overrightarrow{P_3P_4}
+\frac{S_{35}}{r_{35}}\overrightarrow{P_3P_5}
\right)\\
=&0.
\end{aligned}
\end{equation*} 

\item From \eqref{12} we have
\begin{equation}\label{W/r12}
\begin{aligned}
\frac{\partial \tilde W}{\partial r_{12}}
=&S_{12}-64r_{12}(\eta_5\triangle_{234}^{+}\triangle_{134}^{-}+\eta_4\triangle_{235}^{+}\triangle_{135}^{-}+\eta_3\triangle_{245}^{+}\triangle_{145}^{-})\\
=&S_{12}-r_{12}\left(
\frac{S_{34}}{r_{34}}\frac{\triangle_{234}^{+}\triangle_{134}^{-}}{\triangle_{124}^{+}\triangle_{123}^{-}}
+\frac{S_{35}}{r_{35}}\frac{\triangle_{235}^{+}\triangle_{135}^{-}}{\triangle_{125}^{+}\triangle_{123}^{-}}
+\frac{S_{45}}{r_{45}}\frac{\triangle_{245}^{+}\triangle_{145}^{-}}{\triangle_{125}^{+}\triangle_{124}^{-}}\right)\\
=&\frac{r_{12}}{\triangle_{123}^{-}\triangle_{124}^{+}\triangle_{125}^{+}}
\left(
\frac{S_{12}}{r_{12}}\triangle_{123}^{-}\triangle_{124}^{+}\triangle_{125}^{+}
-\frac{S_{34}}{r_{34}}\triangle_{234}^{+}\triangle_{134}^{-}\triangle_{125}^{+}\right.\\
&\left.\qquad\qquad\qquad\quad
-\frac{S_{35}}{r_{35}}\triangle_{235}^{+}\triangle_{135}^{-}\triangle_{124}^{+}
+\frac{S_{45}}{r_{45}}\triangle_{245}^{+}\triangle_{145}^{-}\triangle_{123}^{-}
\right).
\end{aligned}
\end{equation}
Now by using \eqref{x1} we have
$$\frac{S_{12}}{r_{12}}\overrightarrow{P_2P_1}=\frac{S_{13}}{r_{13}}\overrightarrow{P_1P_3}+\frac{S_{14}}{r_{14}}\overrightarrow{P_1P_4}+\frac{S_{15}}{r_{15}}\overrightarrow{P_1P_5}.$$
From \eqref{W/r15} and \eqref{W/r14}, we have
\begin{equation*}
\begin{aligned}
\frac{S_{15}}{r_{15}}\triangle_{125}^{+}
-\frac{S_{35}}{r_{35}}\triangle_{235}^{+}
-\frac{S_{45}}{r_{45}}\triangle_{245}^{+}&=0,\\
\frac{S_{14}}{r_{14}}\triangle_{124}^{+}
-\frac{S_{34}}{r_{34}}\triangle_{234}^{+}
+\frac{S_{45}}{r_{45}}\triangle_{245}^{+}&=0.\\
	\end{aligned}
\end{equation*}
To deal with the first term in the parenthesis on the right-hand side of the last equation in \eqref{W/r12}, we compute
\begin{equation}\label{W/r12second}
\begin{aligned}
&\frac{S_{12}}{r_{12}}\triangle_{123}^{-}\triangle_{124}^{+}\triangle_{125}^{+}
=\triangle_{124}^{+}\triangle_{125}^{+}\frac{S_{12}}{r_{12}}\left(-\overrightarrow{P_1P_2}\times \overrightarrow{P_1P_3}\right)\\
=&\triangle_{124}^{+}\triangle_{125}^{+}
\left(
\frac{S_{13}}{r_{13}}\overrightarrow{P_1P_3}
+\frac{S_{14}}{r_{14}}\overrightarrow{P_1P_4}
+\frac{S_{15}}{r_{15}}\overrightarrow{P_1P_5}
\right)\times \overrightarrow{P_1P_3}\\
=&\triangle_{124}^{+}\triangle_{125}^{+}
\left(
\frac{S_{14}}{r_{14}}\triangle_{143}^{+}
+\frac{S_{15}}{r_{15}}\triangle_{153}^{+}
\right)\\
=&\triangle_{125}^{+}\triangle_{143}^{+}\frac{S_{14}}{r_{14}}\triangle_{124}^{+}
+\triangle_{124}^{+}\triangle_{153}^{+}\frac{S_{15}}{r_{15}}\triangle_{125}^{+}\\
=&\triangle_{125}^{+}\triangle_{143}^{+}
\left(\frac{S_{34}}{r_{34}}\triangle_{234}^{+}
-\frac{S_{45}}{r_{45}}\triangle_{245}^{+}\right)
+\triangle_{124}^{+}\triangle_{153}^{+}
\left(\frac{S_{35}}{r_{35}}\triangle_{235}^{+}
+\frac{S_{45}}{r_{45}}\triangle_{245}^{+}\right)\\
=&\frac{S_{34}}{r_{34}}\triangle_{125}^{+}\triangle_{143}^{+}\triangle_{234}^{+}
+\frac{S_{35}}{r_{35}}\triangle_{124}^{+}\triangle_{153}^{+}\triangle_{235}^{+}
+\frac{S_{45}}{r_{45}}\left(-\triangle_{125}^{+}\triangle_{143}^{+}+\triangle_{124}^{+}\triangle_{153}^{+}\right)\triangle_{245}^{+}.
\end{aligned}
\end{equation}
Substituting the last expression in \eqref{W/r12second} to the last expression in \eqref{W/r12}, we have 
\begin{equation*}
\begin{aligned}
\frac{\partial \tilde W}{\partial r_{12}}
=&\frac{r_{12}}{\triangle_{123}^{-}\triangle_{124}^{+}\triangle_{125}^{+}}\cdot
\frac{S_{45}}{r_{45}}
\left(
-\triangle_{125}^{+}\triangle_{143}^{+}
+\triangle_{124}^{+}\triangle_{153}^{+}
+\triangle_{145}^{-}\triangle_{123}^{-}
\right)\triangle_{245}^{+}.
\end{aligned}
\end{equation*}

Now, we need to show that
\begin{equation*}
-\triangle_{125}^{+}\triangle_{143}^{+}
+\triangle_{124}^{+}\triangle_{153}^{+}
+\triangle_{145}^{-}\triangle_{123}^{-}
=
\triangle_{125}^{+}\triangle_{134}^{+}
-\triangle_{124}^{+}\triangle_{135}^{+}
+\triangle_{145}^{+}\triangle_{123}^{+}=0.
\end{equation*}
Noticing that
\begin{equation*}
	\begin{aligned}
&\triangle_{125}^{+}\triangle_{134}^{+}
-\triangle_{124}^{+}\triangle_{135}^{+}
+\triangle_{145}^{+}\triangle_{123}^{+}\\
=&\left(\overrightarrow{P_1P_2}\times \overrightarrow{P_1P_5}\right)\left(\overrightarrow{P_1P_3}\times \overrightarrow{P_1P_4}\right)
-
\left(\overrightarrow{P_1P_2}\times \overrightarrow{P_1P_4}\right)\left(\overrightarrow{P_1P_3}\times \overrightarrow{P_1P_5}\right)\\
&+
\left(\overrightarrow{P_1P_4}\times \overrightarrow{P_1P_5}\right)\left(\overrightarrow{P_1P_2}\times \overrightarrow{P_1P_3}\right)\\
=&r_{12}r_{13}r_{14}r_{15}
\left(
\sin\varphi_{25}\sin\varphi_{34}
-\sin\varphi_{24}\sin\varphi_{35}
+\sin\varphi_{45}\sin\varphi_{23}
\right),
	\end{aligned}
\end{equation*}
where $\varphi_{ij}\in[0,2\pi)$ is the angle formed by $\overrightarrow{P_1P_i}$ rotating anti-clockwise to $\overrightarrow{P_1P_j}$, and  both of the vectors start from $P_1$. Then there comes an order of the particles scanned anti-clockwise by $\overrightarrow{P_1P_2}$. Without loss of generality, we suppose that the order is $P_2, P_3, P_4, P_5$, and one can see Figure \ref{Wr12.png}. Furthermore, let
$\varphi_{23}=\alpha, \varphi_{34}=\beta, \varphi_{45}=\gamma,$
then 
\begin{equation*}
\begin{aligned}
&\sin\varphi_{25}\sin\varphi_{34}
-\sin\varphi_{24}\sin\varphi_{35}
+\sin\varphi_{45}\sin\varphi_{23}\\
=&\sin(\alpha+\beta+\gamma)\sin\beta-\sin(\alpha+\beta)\sin(\beta+\gamma)+\sin\gamma\sin\alpha\\
=&\sin(\alpha+\beta)\cos\gamma\sin\beta+\cos(\alpha+\beta)\sin\gamma\sin\beta\\
&-\sin(\alpha+\beta)\sin\beta\cos\gamma-\sin(\alpha+\beta)\cos\beta\sin\gamma+\sin\gamma\sin\alpha\\
=&\sin\gamma\sin(\beta-(\alpha+\beta))+\sin\gamma\sin\alpha\\
=&0.
\end{aligned}
\end{equation*}
Hence from \eqref{W/r12} we have $\frac{\partial \tilde W}{\partial r_{12}}=0$.
\end{enumerate}
\begin{figure}[htbp]
	\centering
	\includegraphics[width=0.5\textwidth]{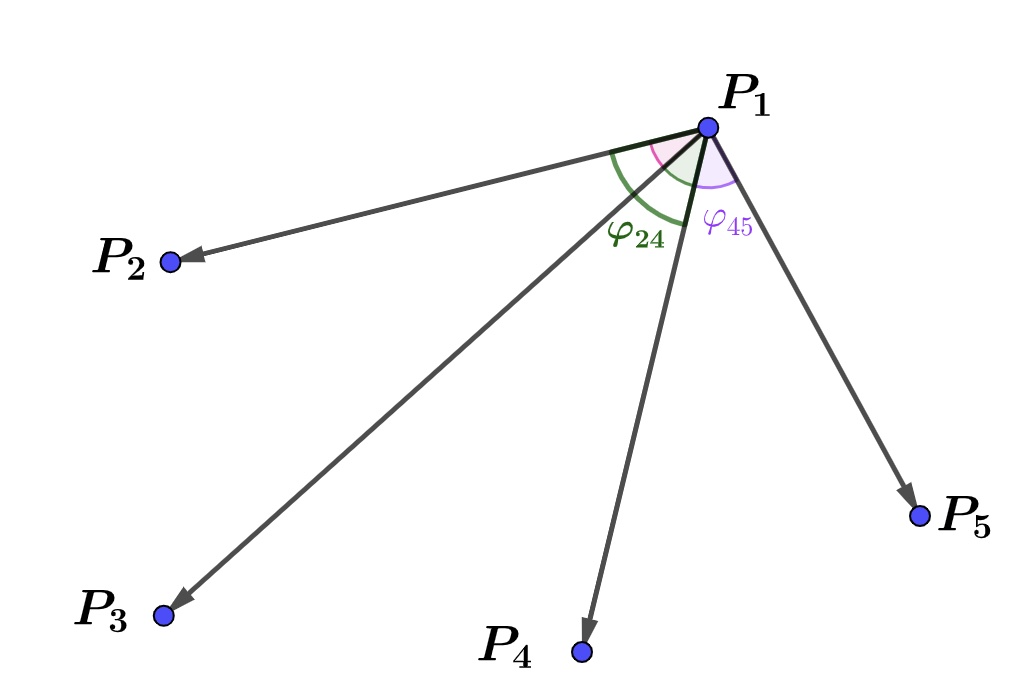}
	\caption{Possible positions of $P_1,\cdots,P_5$.}\label{Wr12.png}
\end{figure}

\begin{rema}
It is well known that except for the 5-body collinear case, there is no 4-body collinearity in the planar 5-body central configurations by the Perpendicular Bisector Theorem in \cite{moeckel1990}. Hence, there are at most two 3-body collinear cases in one planar 5-body central configuration. For example, for the trapezoidal shape, if we choose a constraint $F_i$ involving a 3-body collinearity, we are unable to find probable $\eta_4$ and $\eta_5$ to fit \eqref{12}-\eqref{25} by direct computation. That means $F_i$ is of no avail; in other words, $\tilde W(r)$ is not the equivalent Lagrangian function for this shape. The same happens if there are two 3-body collinear cases in the planar 5-body central configuration, for example, the rhombus shape with one mass in the center. 
\end{rema}

\section*{Acknowledgements}
The authors sincerely thank the anonymous referees for their many valuable comments that helped us improve this paper.
This work is partially supported by the NSF of China (No. 12071316). The authors would like to thank the NSF of China. 

\section*{Data availability}
All data generated or analyzed during this study are included in this published article.

\section*{Conflict of interest statement}

The authors declared that they have no conflicts of interest in this work.

\bibliographystyle{amsplain}
\bibliography{unique.bib}

\end{document}